\RequirePackage{fix-cm}

\documentclass[smallextended]{svjour3}
\smartqed

\usepackage{amsmath}
\usepackage{lipsum}
\usepackage{amsfonts}
\usepackage{graphicx}
\usepackage{epstopdf}
\usepackage{algorithmic}
\usepackage{bm}
\usepackage{booktabs}
\usepackage{xcolor}
\usepackage{url}
\usepackage{orcidlink}

\usepackage{arydshln}

\usepackage{amsopn}

\UseRawInputEncoding
\usepackage[utf8]{inputenc}

\ifpdf
  \DeclareGraphicsExtensions{.eps,.pdf,.png,.jpg}
\else
  \DeclareGraphicsExtensions{.eps}
\fi

\usepackage{enumitem}
\setlist[enumerate]{leftmargin=.5in}
\setlist[itemize]{leftmargin=.5in}

\begin{document}

\title{Fast numerical solvers for parameter identification problems in mathematical biology
}

\titlerunning{Fast solvers for parameter identification problems}

\author{Karol\'{i}na Benkov\'{a}         \and
        John W. Pearson         \and
        Mariya Ptashnyk
}


\institute{Karol\'{i}na Benkov\'{a} \at
              The Maxwell Institute for Mathematical Sciences, The Bayes Centre \\ 47 Potterrow, Edinburgh, EH8 9BT, Scotland, UK \\
              \email{kbenkova@ed.ac.uk}           %  \\
           \and
           John W. Pearson \at
              School of Mathematics, The University of Edinburgh \\ The King’s Buildings, Peter Guthrie Tait Road, Edinburgh, EH9 3FD, Scotland, UK \\
              \email{j.pearson@ed.ac.uk}
           \and
           Mariya Ptashnyk \at
              School of Mathematical and Computer Sciences, Heriot-Watt University \\ Riccarton Campus, Edinburgh, EH14 4AS, Scotland, UK \\
              \email{m.ptashnyk@hw.ac.uk}
}

\date{Received: date / Accepted: date}

\maketitle

\begin{abstract}
In this paper, we consider effective discretization strategies and iterative solvers for nonlinear PDE-constrained optimization models for pattern evolution within biological processes. Upon a Sequential Quadratic Programming linearization of the optimization problem, we devise appropriate time-stepping schemes and discrete approximations of the cost functionals such that the discretization and optimization operations are commutative, a highly desirable property of a discretization of such problems. We formulate the large-scale, coupled linear systems in such a way that efficient preconditioned iterative methods can be applied within a Krylov subspace solver. Numerical experiments demonstrate the viability and efficiency of our approach.
\end{abstract}

\section{Introduction}
\label{sec:intro}

To deepen our understanding of complex biological systems, it is crucial to develop mathematical models to uncover underlying mechanisms and analyze such systems. 
One of the most fundamental challenges in biology is to investigate how structures and shapes of organisms arise \cite{greensharpe2015}, or how an organism develops from, e.g., a single fertilized egg \cite{meinhardt2012}, which involves modelling and studying developmental processes.
A crucial mechanism in embryonic development is the formation of spatial patterns that then evolve into formation of different tissues and organs. 
Modelling the mechanisms that give rise to these patterns is an active research area and a variety of structures have been studied, ranging from patterns on the skin that give rise to hair follicles \cite{glover2017}, to the formation of organs such as the villi of the intestine~\cite{walton2016} or the development of digits of a limb~\cite{raspopovic2014}.

The most well-known theory for spatial pattern formation was first suggested by Alan Turing in his seminal  paper \cite{turing}. 
Turing observed that asymmetric organisms develop from symmetric embryos, and  the symmetry is broken during morphogenesis, the process of tissue and organ formation in an embryo.
He proposed that this ``symmetry-breaking'' is caused by  reactions between the chemicals involved, termed the \emph{morphogens}, as their interactions can amplify small disturbances in the nearly homogeneous tissue. 
By adding diffusion to each morphogen at different rates, instability is added to the system, which gives rise to so-called \emph{Turing patterns}.
While there exist other theories for pattern formation, Turing's reaction--diffusion theory is still a widely used modelling tool \cite{painterptashnyk2022,sudderick2024}.

The dynamics of morphogens  can be modelled using a system of partial differential equations (PDEs). 
To analyse and simulate the system numerically, we first need to determine the parameters and source functions included in the equations.
In an ideal world, all of these can be measured directly,
however, in practice we are frequently required to fit them using indirect measurements. 
A  parameter identification problem consists of solving an inverse problem.
Assuming we have a candidate for one or more parameters (or functions) within a PDE model describing a biological system, along with one or more observations such as images, 
a valuable approach for determining or approximating the missing parameters is to formulate an optimization problem with these PDEs serving as constraints and the unknown parameters as controls. This is an example of a \emph{PDE-constrained optimization} (PDECO) problem; see \cite{MQS,troltzsch} for comprehensive overviews of such problems. The optimization aspect consists of minimizing the deviation of solutions of the PDEs from the observation while also minimizing the amount of control required to do so.

In order to solve such a PDECO problem numerically, we derive a system of discrete optimality conditions in the form of discretized PDEs. 
While there are a variety of approaches for solving PDE systems, a key challenge when tackling time-dependent problems is that the coupled system involves PDEs evolving forward and backward in time. 
For this reason, and in order to achieve good convergence properties, we elect to solve the system with the \emph{all-at-once method}, that is we simultaneously solve the linearized and discretized system of optimality conditions for all variables, across the entire time interval and within the whole spatial domain. 
Clearly, the dimensions of the linearized systems for fine mesh and time interval discretizations can be extremely large, which necessitates a highly efficient and robust solution strategy, with minimal computational storage requirements. We therefore solve the fully discretized system iteratively, and devise a suitable preconditioning strategy to accelerate convergence and allow us to obtain the solution fast and robustly. 
 We highlight that the all-at-once approach for solving reaction–diffusion models, employing backward Euler discretization in time, was demonstrated in \cite{stollpearsonmaini2016}. 
Another class of methods for addressing such problems are gradient-based techniques. 
Notable examples of image-driven PDECO problems for reaction-diffusion equations include \cite{hogea2008}, and \cite{garvie2010,garvie2014,sgura2018} for spatial patterns in particular. Additionally, similar approaches have been applied to reaction--diffusion--advection equations~\cite{uzunca2017}.

The process of transforming the optimization problem into a fully discretized PDE system to be solved numerically can be done by performing the optimization and discretization steps in either order.
Developing \emph{adjoint-consistent} approaches that enable discretization and optimization to commute, thus making the system matrices resulting from the Optimize-then-Discretize (OTD) and Discretize-then-Optimize (DTO) approaches equivalent, is of great interest in the field. It is also of importance when one wishes to quantify the discretization error properties of such problems. 
However, this is a challenging problem even for relatively simple (linear) PDEs, e.g., applying the two approaches to the advection--diffusion equation with streamline upwind/Petrov Galerkin (SUPG) stabilization was shown to lead to different results \cite{collis2002}. The DTO approach results in a symmetric system matrix given appropriate discretization choices, while the symmetry of the OTD matrix depends on the time discretization scheme of choice. Hence, we wish to develop a strategy that enables commutation between the discretization and optimization steps for challenging nonlinear (in particular semilinear) PDE systems arising from modelling pattern formation. Inspired by the St\"{o}rmer--Verlet integrator (see \cite{hairerlubich2003} for a review), we also aim to obtain satisfying convergence properties upon refinement in the time discretization.

The aim of this paper is to propose a numerical scheme to solve PDECO problems with reaction--diffusion equations, which allows the commutation of the discretization and optimization steps.
The discretization scheme that results in a symmetric system matrix also allows us to exploit the robust and theoretically-predictable convergence properties of suitable iterative methods for symmetric systems.
In particular, we wish to derive a fast and efficient iterative numerical solver for the fully discretized problem by exploiting the algebraic structure of the system matrix and designing a suitable preconditioner. 
We highlight that the framework we present 
is general and can be applied to a wide range of nonlinear evolution PDEs. It also opens up the possibility of applying such a strategy to other time-stepping methods for a range of PDEs, which one may have reason to apply in order to achieve desirable stability properties or make use of implicit--explicit (IMEX) methods\footnote{We note that in this latter case, one might sacrifice the adjoint-consistency property at the very start and end of the time interval.}.

The paper is organized as follows. 
In Section \ref{sec:background}, we formulate a parameter identification problem for a system of  reaction--diffusion equations in the PDECO context, 
which is linearized using the Sequential Quadratic Programming method in Section~\ref{sec:nonlinear_programming}. 
The linear--quadratic PDECO problem that emerges is then transformed into a fully discretized system in the following two sections: in Section \ref{sec:otd}, we apply the St\"{o}rmer--Verlet time discretization method to the optimality conditions in the OTD approach, while in Section~\ref{sec:dto} we first approximate the Lagrangian with midpoint and trapezoidal methods and then derive the optimality system. The system is then written in matrix--vector form in Section~\ref{sec:all_at_once}. To be able to solve the system of linear algebraic equations fast, in Section~\ref{sec:precond} we derive a suitable preconditioner. 
In Section~\ref{sec:backward_euler} we introduce the backward Euler method for our PDECO problem and in Section~\ref{sec:results} we compare numerical results for the  St\"{o}rmer--Verlet and classical backward Euler methods. The concluding remarks are presented in Section~\ref{sec:conc}. 

\subsection{Notation}
The matrices are denoted by uppercase letters in bold font.
As the work in this paper involves block diagonal, sub- and super-diagonal matrices, we introduce a concise notation for their representation to enhance brevity in mathematical expressions. We denote block diagonal matrices, with $N$ diagonal blocks that contain matrices $\mathbf{P}^i$ or $\mathbf{R}^i$, by
\begin{equation*}
	\text{blkdiag} \left( \left\lbrace \mathbf{P}^i\right\rbrace_{i=1}^N \right) \hspace{3mm} \text{or}  \hspace{3mm} \text{blkdiag} \left( \left\lbrace \mathbf{P}^i\right\rbrace_{i=1}^{N_1}, \left\lbrace \mathbf{R}^i\right\rbrace_{i=1}^{N_2} \right),
\end{equation*}
where the superscripts indicate the time-step of the elements involved in the block matrix 
and $N_1+N_2=N$. Similarly, we define the operators for block sub- and super-diagonal matrices, with sub-diagonals $\mathbf{U}^i$ and super-diagonals $\mathbf{V}^i$, 
as
\begin{equation*}
	\text{blksubdiag} \left( \left\lbrace \mathbf{U}^i\right\rbrace _{i=1}^{N-1} \right)\hspace{3mm} \text{and} \hspace{3mm} \text{blksupdiag} \left( \left\lbrace \mathbf{V}^i\right\rbrace _{i=1}^{N-1} \right) ,
\end{equation*}
respectively. Calligraphic typeface denotes the main system matrices and their preconditioners.

The vectors are represented by lowercase letters in bold font. To differentiate between various types of vectors utilized in different contexts, continuous variable vectors and vector functions are represented using bold roman typeface, e.g., $\textbf{u} = [u,v]^\top$, whereas vectors containing fully discretized variables are indicated by bold italic, e.g., $\bm{u}^i$ represents the variable $u$ discretized over the mesh at the $i$-th time-step. Underlined bold italic vectors combine two variables discretized over the mesh and across some time interval, e.g., $\underline{\bm{u}} := \left[(\bm{u}^1)^{\top}, ..., (\bm{u}^{N_t})^{\top}, (\bm{v}^1)^{\top}, ..., (\bm{v}^{N_t})^{\top}\right]^{\top}$ for time-steps $1, ..., N_t$.

\section{Parameter identification for reaction--diffusion equations}
\label{sec:background}

To determine or approximately identify missing parameters or functions in a system of PDEs, we may formulate a constrained optimization problem with the PDEs as constraints. 
For the problems considered here, we think of the missing parameters as \emph{control variables}, and aim to find optimal controls which minimize the misfit between the solutions of the PDEs (\emph{state variables}) and observed data in the $L^2$-norm, which is a frequent choice of norm for these problems.
To ensure the existence of a minimum in the optimal control problem, we add to the cost functional defining the optimization problem a regularization term, minimizing the $L^2$-norm of controls.
We are interested in finding solutions in a bounded Lipschitz domain $\Omega \subset \mathbb{R}^2$ throughout the time interval $(0,T)$ and  we solve the optimization problem in the time--space cylinder $Q:=(0,T)\times\Omega $. The boundary of the spatial domain is denoted by $\partial \Omega$ and the time--space boundary is $\Gamma:=  (0,T)\times \partial\Omega $.

For a PDECO problem with two state variables $u$, $v$ and controls $a$, $b$, we consider the cost functional $\mathcal{J} = \mathcal{J}(u,v,a,b)$:
\begin{equation*}
\begin{aligned}
    % \mathcal{J}(u,v,a,b)
    \mathcal{J} :={}& \frac{\alpha_1}{2} \|u-\hat{u}\|^2_{L^2(Q)} +\frac{\alpha_2}{2} \|v-\hat{v}\|^2_{L^2(Q)} + \frac{\beta_1}{2} \|a\|^2_{L^2(Q)} + \frac{\beta_2}{2} \|b\|^2_{L^2{(Q)}},
    \label{eq:costfun}
\end{aligned}
\end{equation*}
where  $\hat{u}$, $\hat{v}$ represent the observations, known as desired states. 
For the weights $\alpha_i$ and $\beta_i$ ($i=1,2$) we typically choose $\alpha_i>\beta_i$ to prioritize minimizing the misfit. 
In the case of large~$\beta_i$, one would be penalizing the cost of implementing the controls, which may result in the state variables failing to match the desired states closely.

We formulate the PDECO problem for reaction--diffusion systems 
\begin{equation}
    \min_{u,v,a,b} \mathcal{J}(u,v,a,b)
    \label{eq:schnakenberg_min}
\end{equation}
subject to:
\begin{equation}
    \begin{aligned}
        u_t - D_u \nabla^2u + \Phi(u,v) &= \gamma a  \hspace{5mm} \text{ in } Q, \\
        v_t - D_v \nabla^2v + \Psi(u,v) &= \gamma b \hspace{5.3mm} \text{ in } Q,  
    \end{aligned}
     \label{eq:schnakenberg_eq}
\end{equation}
with initial and boundary conditions
\begin{equation}
	\begin{aligned}
		u(0,\textbf{x}) = u_0(\textbf{x}), \hspace{2mm}  v(0,\textbf{x}) = v_0(\textbf{x})\hspace{5mm} &\text{ in } \Omega, \\
		\frac{\partial u}{\partial \bm{\nu}} = \frac{\partial v}{\partial \bm{\nu}} = 0 \hspace{11.4mm} &\text{ on } \Gamma,
	\end{aligned}
	\label{eq:schnakenberg_BC}
\end{equation}
where $\Phi(u,v)$ and $\Psi(u,v)$ are nonlinear reaction terms and $\gamma$ is a positive constant.
The parameters $D_u$ and $D_v$ denote the positive diffusion coefficients and, for the occurrence of Turing patterns, we require $D_u<D_v$. The given positive parameter $\gamma$ can be interpreted as the relative strength of the reaction terms, which for the problem we consider is also a factor within $\Phi$ and $\Psi$, while the control variables $a$ and $b$ represent the production rates for $u$ and $v$, respectively. In this paper, we only consider the case where the control variables are acting on the whole domain, i.e., the controls are distributed, and where they act as source terms in the model equations.
The initial conditions $u_0$ and $v_0$ are assumed to be known.
While we consider zero-flux boundary conditions with the vector normal to the boundary denoted by $\bm{\nu}$, it is also possible to use Dirichlet boundary conditions instead. 

The numerical scheme that we present in this paper is built upon the Schnakenberg model, a system of reaction--diffusion equations commonly used to model pattern formation, for which  $\Phi(u,v)= \gamma(u-u^2v)$ and $\Psi(u,v)= \gamma u^2v$.
The equations were initially proposed by Schnakenberg to study limit cycles arising in chemical reactions \cite{schnakenberg1979}. They have since also been used to study self-organization phenomena in various biological systems, for example, the emergence of skin patterns \cite{shaw1990} or initiation of plant root hair \cite{brenamedina2014}. 

Since the reaction terms are nonlinear, the optimization problem \eqref{eq:schnakenberg_min}--\eqref{eq:schnakenberg_BC} is also nonlinear. To resolve this, we need to utilize some nonlinear programming algorithm. In the following section, we derive the first-order optimality conditions using the Lagrange method and apply the Sequential Quadratic Programming (SQP) method; see \cite[Sec. 5.9.2]{troltzsch} and \cite[Sec. 3.2]{HK10}.

\section{Nonlinear programming}
\label{sec:nonlinear_programming}
Below, we perform linearization for the nonlinear problem \eqref{eq:schnakenberg_min}--\eqref{eq:schnakenberg_BC} with the SQP method, making use of a Newton linearization for the PDE constraints. We commence by following the OTD approach, and via our specific choice of time-stepping and approximation of certain terms, we find that the same optimality conditions are obtained from both OTD and DTO formulations, thus achieving that the discretization and optimization operations are commuting.

\subsection{Lagrange method}
\label{sec:lagrange_method}
Let us start by formulating the continuous Lagrangian $\mathcal{L} = \mathcal{L}(u,v,a,b,p,q)$. We multiply the PDE constraints by the Lagrange multipliers  and add them to the cost functional $\mathcal{J}$:
\begin{equation}
	\begin{split}
		\mathcal{L} (u,v,a,b,p,q):=  \mathcal{J} + \int_Q p \left (u_t - D_u \nabla^2u + \Phi(u,v) - \gamma a \right ) dQ 
		\\
  + \int_Q q \left (v_t - D_v \nabla^2v + \Psi(u,v) - \gamma b \right )  dQ,
	\end{split}
	\label{eq:OTD_Lagrangian}
\end{equation}
where $p$ and $q$ denote the multipliers, also called the \emph{adjoint variables}. We omit including the boundary conditions within the Lagrangian \eqref{eq:OTD_Lagrangian} and the following sections, for the sake of brevity of presentation. The following optimality conditions are a result of finding the stationary points of $\mathcal{L}$, which we obtain by taking Fréchet derivatives; for a comprehensive description of this process for a range of problems we refer to \cite{troltzsch}. 

Firstly, from the partial derivatives with respect to the state variables $u$ and $v$, we derive the adjoint equations:
\begin{equation*}
    \begin{aligned}
    -p_t - D_u \nabla^2 p + \Phi_u(u,v) p  + \Psi_u(u,v) q + \alpha_1(u-\hat{u}) &= 0 \hspace{5mm} \text{in } Q, \\
        -q_t - D_v \nabla^2 q + \Phi_v(u,v) p  + \Psi_v(u,v) q + \alpha_2(v-\hat{v})&= 0 \hspace{5mm} \text{in } Q, 
    \end{aligned}
    %\label{eq:OTD_adjoint}
\end{equation*}
equipped with final-time conditions and zero-flux boundary conditions:
\begin{equation}
    \begin{aligned}
         p(T,\textbf{x}) =  q(T,\textbf{x}) &= 0 \hspace{5mm} \text{in } \Omega, \\
         \frac{\partial p}{\partial \bm{\nu}} = \frac{\partial q}{\partial \bm{\nu}} &= 0 \hspace{5mm} \text{on } \Gamma.
    \end{aligned}
    \label{eq:OTD_BC}
\end{equation}
Next, the derivatives with respect to the controls give us the gradient or control equations:
\begin{equation}
    \begin{aligned}
         \beta_1 a - \gamma p &= 0 \hspace{5mm} \text{in } Q,\\
         \beta_2 b - \gamma q &= 0 \hspace{5mm} \text{in } Q.
    \end{aligned}
    \label{eq:OTD_gradient}
\end{equation}
We note that the gradient equations are algebraic for this problem, and thus can be used to reduce the complexity of the optimization problem in later sections by substituting for the variables $a$ and $b$.
Finally, the set of optimality conditions is completed by studying optimality with respect to the Lagrange multipliers: the Fréchet derivatives will yield the original state equations and initial and boundary conditions as in \eqref{eq:schnakenberg_eq}--\eqref{eq:schnakenberg_BC}.

\subsection{Lagrange--SQP method}
\label{sec:lagrange_sqp}
After obtaining the optimality conditions, we proceed to their linearization.
For our optimality system,
assuming previously calculated Newton iterates at the $k$-th iteration, marked in the subscripts of the variables, we need to solve for the variables $u$, $v$, $p$, $q$ at the current $(k+1)$-th Newton iteration.
Using the more concise notation $\bar \Phi := \Phi(u_k,v_k)$, $\bar \Psi := \Psi(u_k,v_k)$, $\bar \Phi_{u} := \partial_{u}\Phi(u_k,v_k)$, $\bar \Psi_{uv} := \partial^2_{uv} \Psi(u_k,v_k)$, and similarly for other derivatives (where the subscripts of $\Phi$, $\Psi$ denote partial derivatives), we thus need to solve the linearized state equations:
\begin{equation}
	\begin{aligned}
		u_t - D_u \nabla^2 u  + \bar \Phi + \bar \Phi_u(u-u_k) + \bar \Phi_v(v-v_k) - \frac{\gamma^2}{\beta_1} p &= 0, \\
		v_t - D_v\nabla^2v + \bar \Psi + \bar \Psi_u(u-u_k) + \bar \Psi_v(v-v_k)- \frac{\gamma^2}{\beta_2} q & = 0, \\ 
	\end{aligned}
		\label{eq:OTD_Newt_st}
\end{equation}
and, upon an application of the SQP method, the linearized adjoint equations of the form:
\begin{equation}
	\begin{aligned}
		&-p_t- D_u\nabla^2 p + \bar \Phi_u p + \bar \Psi_u q + (\alpha_1 + \bar \Phi_{uu} p_k + \bar\Psi_{uu} q_k) u + (\bar \Phi_{uv} p_k + \bar \Psi_{uv} q_k) v \\ 
		&\hspace{5mm} = \alpha_1 \hat{u} + (\bar \Phi_{uu} p_k + \bar \Psi_{uu} q_k)u_k + (\bar \Phi_{uv} p_k + \bar \Psi_{uv} q_k) v_k, \\ 
		&- q_t - D_v \nabla^2 q + \bar \Phi_v p + \bar \Psi_v q + (\bar \Phi_{vu} p_k + \bar \Psi_{vu} q_k) u + (\alpha_2 + \bar \Phi_{vv} p_k + \bar \Psi_{vv} q_k) v \\
		&\hspace{5mm} = \alpha_2 \hat{v} + (\bar \Phi_{vu} p_k + \bar \Psi_{vu} q_k) u_k + (\bar \Phi_{vv} p_k + \bar \Psi_{vv} q_k) v_k.
	\end{aligned}
	\label{eq:OTD_Newt_ad}  
\end{equation}
The linearized equations \eqref{eq:OTD_Newt_st}--\eqref{eq:OTD_Newt_ad} are equivalent to the first-order optimality conditions of a linear--quadratic optimization problem, obtained via an iteration of the SQP method, as described in \cite{troltzsch}. 
The SQP cost functional $\tilde{\mathcal{J}}$ is related to the cost functional and Lagrangian of the original problem \eqref{eq:schnakenberg_min}--\eqref{eq:schnakenberg_BC} as follows:
\begin{equation}
	\begin{aligned}
		\tilde{\mathcal{J}}(u,v,a,b):= \mathcal{J}'(u_k,v_k,a_k,b_k)\textbf{w} + \frac{1}{2} \textbf{w} ^\top \mathcal{L}''(u_k,v_k,a_k,b_k,p_k,q_k)\textbf{w},
	\end{aligned}
	\label{eq:SQP_Jtilde}
\end{equation}
where  $\textbf{w}=[u-u_k, v-v_k, a-a_k,b-b_k]^\top$ and $\mathcal{L}''$ is the Hessian of~\eqref{eq:OTD_Lagrangian} with respect to $u$, $v$, $a$, and $b$.
At each SQP iteration, we therefore solve the quadratic optimization problem
\begin{equation}
	\min_{u,v,a,b} \mathcal{\tilde{J}}(u,v,a,b) \\
	\label{eq:SQP_min}
\end{equation}
subject to the system of linear reaction--diffusion equations:
\begin{equation}
    \begin{aligned}
     u_t - D_u \nabla^2 u  + \bar \Phi + \bar \Phi_u(u-u_k) + \bar \Phi_v(v-v_k) - \gamma  a &= 0  \hspace{11.4mm} \text{ in } Q, \\
     v_t - D_v\nabla^2v + \bar \Psi + \bar \Psi_u(u-u_k) + \bar \Psi_v(v-v_k)- \gamma  b &= 0  \hspace{11.4mm} \text{ in } Q, \\
     u(0,\textbf{x}) = u_0(\textbf{x}), \hspace{2mm}  v(0,\textbf{x}) &= v_0(\textbf{x})\hspace{5mm} \text{ in } \Omega, \\
     \frac{\partial u}{\partial \bm{\nu}} = \frac{\partial v}{\partial \bm{\nu}} &= 0 \hspace{11.4mm} \text{ on } \Gamma.
    \end{aligned}
     \label{eq:SQP_constraints}
\end{equation}
We now evaluate all the terms  in \eqref{eq:SQP_Jtilde}. The first term is
\begin{equation}
    \begin{aligned}
       \mathcal{J}^\prime(u_k,v_k,a_k,b_k)\textbf{w} ={}&  \alpha_1 \int_Q (u_k-\hat{u})(u-u_k)  \, dQ 
        + \alpha_2 \int_Q (v_k-\hat{v})(v-v_k) \, dQ \\
       & + \beta_1 \int_Q a_k(a-a_k) \, dQ + \beta_2 \int_Q  b_k(b-b_k)  \, dQ,
    \end{aligned}
    \label{eq:SQP_JI}    
\end{equation}
and the second term is
\begin{equation}
    \begin{aligned}
    \frac{1}{2}
     \textbf{w}^\top \mathcal{L}^{\prime \prime}(u_k,v_k,a_k,b_k)\textbf{w} =
        \frac{1}{2} 
        \int_Q
        \textbf{w}^\top \mathbf{D}_{\alpha, \beta} \textbf{w} \, dQ + \frac{1}{2}
        \int_Q
        \textbf{w}^\top \mathbf{A}_{\Phi,\Psi}  \textbf{w} \, dQ,
    \end{aligned}
    \label{eq:SQP_LII}
\end{equation}
where 
$$
\mathbf{D}_{\alpha, \beta}  = \begin{bmatrix} 
   \alpha_1 \text{id} & 0 & 0 & 0 \\
            0 & \alpha_2 \text{id} & 0 & 0 \\
            0 & 0 & \beta_1 \text{id} & 0 \\
            0 & 0 & 0 & \beta_2 \text{id}
\end{bmatrix}, 
~ % \quad%!! 
\mathbf{A}_{\Phi,\Psi}  = \begin{bmatrix} 
p_k \bar \Phi_{uu} + q_k \bar \Psi_{uu} & p_k \bar \Phi_{uv} + q_k \bar \Psi_{uv} & 0 & 0 \\
	p_k \bar \Phi_{vu} + q_k \bar \Psi_{vu} & p_k \bar \Phi_{vv}+ q_k \bar \Psi_{vv}& 0  & 0  \\
	0 & 0 & 0 & 0 \\
	0 & 0 & 0 &  0
\end{bmatrix}, 
$$
and $\text{id}$ denotes the identity operator.

In the remainder of the work, we consider reaction--diffusion equations specific to the Schnakenberg model, but the approach considered here can be applied to any system in the  form~\eqref{eq:schnakenberg_eq}. 
The cost functional \eqref{eq:SQP_Jtilde} can be split into two terms:  
\begin{equation*}
    \mathcal{\tilde{J}}(u,v,a,b) = \int_Q  (\mathcal{\tilde{J}}_1 + \mathcal{\tilde{J}}_2) \, dQ, 
\end{equation*}
where $\mathcal{\tilde{J}}_1$ contains the terms from  \eqref{eq:SQP_JI} and the first integral in \eqref{eq:SQP_LII}:
\begin{equation*}
    \begin{aligned}
        \mathcal{\tilde{J}}_1 ={}& \alpha_1 (u_k-\hat{u})(u-u_k) + \alpha_2 (v_k-\hat{v})(v-v_k) +  \beta_1 a_k(a-a_k) + \beta_2 b_k(b-b_k) \\
        &+ \frac{\alpha_1}{2}(u-u_k)^2 + \frac{\alpha_2}{2}(v-v_k)^2 + \frac{\beta_1}{2}(a-a_k)^2 + \frac{\beta_2}{2}(b-b_k)^2
    \end{aligned}
\end{equation*}
and 
$\mathcal{\tilde{J}}_2$ denotes the second term in \eqref{eq:SQP_LII}:
\begin{equation*}
        \mathcal{\tilde{J}}_2 = \frac{1}{2} \textbf{w}^{\top} \mathbf{A}_{\Phi,\Psi} \textbf{w}, \quad 
        \text{ where } \; \; 
        \mathbf{A}_{\Phi,\Psi} =  \begin{bmatrix}
            2\gamma (q_k-p_k) v_k  & 2\gamma (q_k-p_k) u_k & 0 & 0 \\
            2\gamma (q_k-p_k) u_k & 0 & 0 & 0 \\
            0 & 0 & 0 & 0 \\
            0 & 0 & 0 &  0
            \end{bmatrix}.
\end{equation*}
As $\mathcal{\tilde{J}}$ is a cost functional, and it will later be part of a Lagrangian corresponding to the problem \eqref{eq:SQP_min}--\eqref{eq:SQP_constraints}, the terms independent of $u$, $v$, $a$, $b$ within $\mathcal{\tilde{J}}_1$ and $\mathcal{\tilde{J}}_2$ will not appear in the optimality conditions. We can therefore conveniently add the constant terms $\frac{\alpha_1}{2}\hat{u}^2$ and $\frac{\alpha_2}{2}\hat{v}^2$ to $\mathcal{\tilde{J}}_1$ which allows us to reconstruct the original cost functional $\mathcal{\tilde{J}}$, while the other terms in $\mathcal{\tilde{J}}_1$ will be constant. Specifically,
\begin{equation*}
		\mathcal{\tilde{J}}_1 = \mathcal{J} - \frac{\alpha_1}{2}\hat{u}^2 - \frac{\alpha_2}{2}\hat{v}^2  - \frac{\alpha_1}{2} u_k^2 + \alpha_1 \hat{u}u_k - \frac{\alpha_2}{2} v_k^2 + \alpha_2 \hat{v}v_k -
		\frac{\beta_1}{2}a_k^2 - \frac{\beta_2}{2}b_k^2.
\end{equation*}
Omitting the constant terms, the SQP method consists of minimizing at each iteration the cost functional
\begin{equation}
	\begin{aligned}
		 \mathcal{\tilde{J}} ={}& \frac{\alpha_1}{2} \int_Q (u-\hat{u})^2 \, dQ + \frac{\alpha_2}{2} \int_Q (v-\hat{v})^2 \, dQ + \frac{\beta_1}{2} \int_Q a^2 \, dQ + \frac{\beta_2}{2} \int_Q b^2 \, dQ \\
		&+\int_Q \gamma (q_k-p_k)v_k(u-u_k)^2 \, dQ + 
		\int_Q 2\gamma (q_k-p_k)u_k(u-u_k)(v-v_k) \, dQ,
	\end{aligned}
	\label{eq:SQP_Jtilde_final}
\end{equation}
and we can rearrange the linearized state equations to obtain
\begin{equation}
    \begin{aligned}
     u_t - D_u \nabla^2u + \gamma(1-2u_kv_k)u - \gamma u_k^2 v - \gamma a + 2\gamma u_k^2v_k &= 0 \hspace{5mm} \text{in } Q, \\
      v_t - D_v \nabla^2v + 2\gamma u_kv_k u + \gamma u_k^2v - \gamma b - 2\gamma u_k^2v_k &= 0 \hspace{5mm} \text{in } Q, 
    \end{aligned}
     \label{eq:SQP_state}
\end{equation}
with the initial and boundary conditions~\eqref{eq:schnakenberg_BC}. 
We therefore reformulated the  nonlinear optimization problem \eqref{eq:schnakenberg_min}--\eqref{eq:schnakenberg_BC} as the quadratic optimization problem \eqref{eq:SQP_Jtilde_final}, \eqref{eq:SQP_state}, \eqref{eq:schnakenberg_BC}.
Next is to provide a method for which the discretization and optimization steps commute for the SQP problem.

\section{Optimize-then-Discretize approach}
\label{sec:otd}
To solve the SQP problem~\eqref{eq:SQP_Jtilde_final}, \eqref{eq:SQP_state}, \eqref{eq:schnakenberg_BC} using the Optimize-then-Discretize approach, we apply the Lagrange multiplier technique described in Section~\ref{sec:lagrange_sqp} to derive the optimality conditions. 
With the Schnakenberg equations as constraints, the first two optimality conditions recover the linearized state equations~\eqref{eq:SQP_state}, and the other two are the adjoint equations:
\begin{equation}
	\begin{aligned}
		-p_t - D_u\nabla^2 p + \alpha_1 u
		+ 2\gamma v_k \left(q_k-p_k\right) u + 
		2\gamma u_k \left(q_k-p_k\right) v  
	     + \gamma(1- 2 u_k v_k) p& \\
	    {}+2\gamma u_k v_k q = \alpha_1 \hat{u} + 4\gamma u_k v_k\left(q_k - p_k \right){}&, \\
  - q_t - D_v \nabla^2 q 	
		+ 2\gamma u_k \left(q_k-p_k\right)u + \alpha_2 v
		+ \gamma u_k^2 \left( q - p\right) = \alpha_2 \hat{v} + 
		 2\gamma u_k^2 \left(q_k - p_k\right){}&,
	\end{aligned}
	\label{eq:OTD_ad}
\end{equation}
the same as in \eqref{eq:OTD_Newt_ad}. We also recover the same linear gradient equations as in \eqref{eq:OTD_gradient}, and initial and boundary conditions in \eqref{eq:schnakenberg_BC} and \eqref{eq:OTD_BC}, respectively.

\subsection{St\"{o}rmer--Verlet scheme}
\label{sec:otd_stormer_verlet}
St\"{o}rmer--Verlet type schemes for time-dependent optimization problems with linear PDEs as constraints, in particular the heat equation, were considered in~\cite{apelflaig2012,flaig2013}; see also \cite{meidnervexler2011,koster2011}.  We now introduce this scheme in the context of the nonlinear problems of this paper, denoting
$\textbf{u} = [u,v]^\top$ and $\textbf{p} = [p,q]^\top$ and considering
the state and adjoint equations in the form 
\begin{equation}
	\textbf{u}_t=\textbf{f}(\textbf{u},\textbf{p}), \hspace{5mm} \textbf{p}_t=\textbf{g}(\textbf{u},\textbf{p}).
	\label{eq:sv_pdes}
\end{equation}
Notice that using the gradient equations \eqref{eq:OTD_gradient}, the right-hand side functions in both state and adjoint equations can be written as functions of  $\textbf{u}$ and $\textbf{p}$. We therefore have 
$$
\textbf{f}(\textbf{u},\textbf{p}) = 
	\begin{bmatrix}
f_1(\textbf{u},\textbf{p})\\
f_2(\textbf{u},\textbf{p})
 \end{bmatrix} \quad \text{ and } \quad \textbf{g}(\textbf{u},\textbf{p}) = 
	\begin{bmatrix}
g_1(\textbf{u},\textbf{p})\\
g_2(\textbf{u},\textbf{p})
 \end{bmatrix},
$$
where 
\begin{equation}\label{formul_f_g}
\begin{aligned} 
f_1(\textbf{u},\textbf{p}) ={} &  D_u \nabla^2u - \gamma(1-2u_k v_k)u + \gamma u_k^2 v + \frac{\gamma^2}{\beta_1} p - 2\gamma u_k^2v_k, \\
f_2(\textbf{u},\textbf{p}) ={} & D_v \nabla^2v - 2\gamma u_k v_k u - \gamma u_k^2v + \frac{\gamma^2}{\beta_2} q + 2\gamma u_k^2 v_k, \\
g_1(\textbf{u},\textbf{p}) ={} & -D_u \nabla^2p + \gamma(1-2u _k v_k)p + 2\gamma u_k v_k q + 2\gamma v_k(q_k-p_k)u \\
           & + 2\gamma u_k(q_k-p_k)v 
		   - \alpha_1(\hat{u}-u) - 4\gamma u_kv_k(q_k-p_k),\\
 g_2(\textbf{u},\textbf{p}) ={} & - D_v \nabla^2q - \gamma u_k^2 p + \gamma u_k^2 q + 2\gamma u_k(q_k-p_k) u - \alpha_2(\hat{v}-v) \\
        & - 2\gamma u_k^2(q_k-p_k).
\end{aligned} 
\end{equation}
Similarly to \cite[Eq.~(1.24)]{hairerlubich2003}, we  discretize equations \eqref{eq:sv_pdes} using the time-stepping scheme described as follows:
\begin{equation*}
	\begin{aligned}
		\textbf{p}^{n+\frac{1}{2}} &= \textbf{p}^n + \frac{\tau}{2}\textbf{g}(\textbf{u}^n,\textbf{p}^{n+\frac{1}{2}}), \\
		\textbf{u}^{n+1} &= \textbf{u}^n + \frac{\tau}{2}\left( \textbf{f}(\textbf{u}^n,\textbf{p}^{n+\frac{1}{2}}) + \textbf{f}(\textbf{u}^{n+1},\textbf{p}^{n+\frac{1}{2}})\right),  \\
		\textbf{p}^{n+1} &=  \textbf{p}^{n+\frac{1}{2}} + \frac{\tau}{2}\textbf{g}(\textbf{u}^{n+1},\textbf{p}^{n+\frac{1}{2}}), 
	\end{aligned}
\end{equation*}
for $n=0,1,...,N_t-1$, where $N_t$ is the number of time-steps of size $\tau$, and with superscripts denoting the time indices. By combining the two half time-steps in $\textbf{p}$, we obtain the scheme
\begin{equation}
	\begin{aligned}
		\textbf{p}^{n+\frac{1}{2}} &=  \textbf{p}^{n-\frac{1}{2}} + \frac{\tau}{2}\left( \textbf{g}(\textbf{u}^n,\textbf{p}^{n-\frac{1}{2}}) + \textbf{g}(\textbf{u}^n,\textbf{p}^{n+\frac{1}{2}})\right) , & \text{ for } \; n=1,2,...,N_t-1,\\
			\textbf{u}^{n+1} &= \textbf{u}^n + \frac{\tau}{2}\left( \textbf{f}(\textbf{u}^n,\textbf{p}^{n+\frac{1}{2}}) + \textbf{f}(\textbf{u}^{n+1},\textbf{p}^{n+\frac{1}{2}})\right) , & \text{ for } \; n=0,1,...,N_t-1, \\
	\end{aligned}
	\label{eq:stormer_verlet2}
\end{equation}
with the values at $n=0$ and $n=N_t - \frac{1}{2}$ in $\textbf{p}$ calculated using
\begin{equation}
	\begin{aligned}
		\textbf{p}^{\frac{1}{2}} &=  \textbf{p}^0 + \frac{\tau}{2}\textbf{g}(\textbf{u}^0,\textbf{p}^{\frac{1}{2}}), \\
		\textbf{p}^{N_t} &= \textbf{p}^{N_t-\frac{1}{2}} + \frac{\tau}{2}\textbf{g}(\textbf{u}^{N_t},\textbf{p}^{N_t-\frac{1}{2}}).
	\end{aligned}
	\label{eq:stormer_verletp}
\end{equation}
We note that the expressions for $\textbf{u}^{n+1}$ and $\textbf{p}^{n+\frac{1}{2}}$ in \eqref{eq:stormer_verlet2} resemble trapezoidal rule type schemes, and this fact is used in Section~\ref{sec:dto} when constructing the scheme using the DTO approach.
Due to the final-time condition for the adjoint variables, the values of $\textbf{p}^{N_t}$ in \eqref{eq:stormer_verletp} will be equal to zero.
The same time-stepping will be followed when $g$ and $f$ also involve the previous iterates of the variables, e.g., $\textbf{g}(\textbf{u}^n,\textbf{p}^{n-\frac{1}{2}}) = \textbf{g}(\textbf{u}^n,\textbf{u}^n_k,\textbf{p}^{n-\frac{1}{2}}, \textbf{p}^{n-\frac{1}{2}}_k)$.

Applying the scheme over the time interval $(0,T)$  yields the semi-dis\-cretized state equations  
\begin{equation}
	\begin{aligned}
		&-u^n - \frac{\tau D_u}{2}  \nabla^2 u^n + \frac{\tau \gamma}{2} (1-2u_k^n v_k^n) u^n 
		+ u^{n+1} - \frac{\tau  D_u}{2} \nabla^2 u^{n+1} \\
            &{}+  \frac{\tau \gamma}{2}( 1-2u_k^{n+1} v_k^{n+1}) u^{n+1} 
		- \frac{\tau \gamma}{2} ( u_k^n) ^2 v^n - 
		\frac{\tau \gamma}{2} ( u_k^{n+1}) ^2 v^{n+1} \\ 
            &{}-\tau\gamma a^{n+\frac{1}{2}} = 
	    -\tau\gamma \left[ \left(u_k^n\right)^2 v_k^n + ( u_k^{n+1}) ^2 v_k^{n+1} \right]
	\end{aligned}
	\label{eq:SQP_stateu_semi}
\end{equation}
and
\begin{equation}
	\begin{aligned}
		&\tau\gamma u_k^n v_k^n u^n +
		\tau\gamma u_k^{n+1} v_k^{n+1} u^{n+1}
		-v^n - \frac{\tau  D_v }{2}\nabla^2 v^n + \frac{\tau \gamma}{2}  \left(u_k^n\right)^2 v^n 
		+ v^{n+1} \\
            &{}-  \frac{\tau  D_v}{2} \nabla^2 v^{n+1} 
		+ \frac{\tau \gamma}{2} ( u_k^{n+1}) ^2 v^{n+1} 
		- \tau  \gamma b^{n+\frac{1}{2}} = 
		\tau\gamma \left[ \left(u_k^n\right)^2 v_k^n + \left(u_k^{n+1}\right)^2 v_k^{n+1}  \right],
	\end{aligned}
	\label{eq:SQP_statev_semi}
\end{equation}
for $n=0,1,...,N_t-1$, 
with the initial conditions $u^0 = u(0,\textbf{x})$, $v^0 = v(0,\textbf{x})$ in $\Omega$.
For the semi-discretized adjoint equations,
starting with the first half step, we obtain 
\begin{equation}
	\begin{aligned}
		& \frac{\tau \alpha_1}{2} u^0 
         + \tau \gamma  (q_k^{\frac{1}{2}} - p_k^{\frac{1}{2}}) v_k^0 u^0 
         + \tau \gamma (q_k^{\frac{1}{2}} - p_k^{\frac{1}{2}}) u_k^0 v^0 
         + p^0 - p^{\frac{1}{2}}  - \frac{\tau D_u}{2}\nabla^2 p^{\frac{1}{2}} \\
         &{} + \frac{\tau \gamma}{2}  \left(1- 2 u_k^0 v_k^0\right)  p^{\frac{1}{2}}  + \tau \gamma u_k^0 v_k^0 q^{\frac{1}{2}} = \frac{\tau \alpha_1}{2} \hat{u}^0 + 2 \tau \gamma u_k^0 v_k^0 (q_k^{\frac{1}{2}} - p_k^{\frac{1}{2}})
	\end{aligned}
	\label{eq:SQP_adjointp_semi_p0}
\end{equation}
and
\begin{equation}
	\begin{aligned}
		& \tau \gamma  (q_k^{\frac{1}{2}} - p_k^{\frac{1}{2}})  u_k^0 u^0 
		+ \frac{\tau \alpha_2}{2} v^0 
		  -\frac{\tau \gamma }{2} ( u_k^0 )^2 p^{\frac{1}{2}}   
            + q^0 - q^{\frac{1}{2}} - \frac{\tau D_v}{2}\nabla^2 q^{\frac{1}{2}} \\
            &{}+\frac{\tau \gamma }{2} \left(u_k^0 \right)^2  q^{\frac{1}{2}} =  \frac{\tau \alpha_2}{2} \hat{v}^0 
		+ \tau \gamma \left(u_k^0 \right)^2
		(q_k^{\frac{1}{2}} - p_k^{\frac{1}{2}}).
	\end{aligned}
	\label{eq:SQP_adjointq_semi_q0}
\end{equation}
The adjoint equations for $n=1,...,N_t-1$ are
\begin{equation}
	\begin{aligned}
		&\tau \alpha_1 u^n + \tau \gamma \left [ (q_k^{n-\frac{1}{2}} - p_k^{n-\frac{1}{2}}) +
		 (q_k^{n+\frac{1}{2}} - p_k^{n+\frac{1}{2}}) \right ] v_k^n u^n + \tau \gamma \left[ (q_k^{n-\frac{1}{2}} - p_k^{n-\frac{1}{2}}) \right.\\
		 &\left. {}+ (q_k^{n+\frac{1}{2}} - p_k^{n+\frac{1}{2}}) \right]  u_k^n v^n 
		+ p^{n-\frac{1}{2}} - \frac{\tau D_u}{2}\nabla^2 p^{n-\frac{1}{2}} + \frac{\tau \gamma }{2} \left(1- 2 u_k^n v_k^n\right)  p^{n-\frac{1}{2}} \\
            &{}- p^{n+\frac{1}{2}} 
		  - \frac{\tau D_u}{2}\nabla^2 p^{n+\frac{1}{2}} + \frac{\tau \gamma}{2}  \left(1- 2 u_k^n v_k^n\right)  p^{n+\frac{1}{2}} 
		+ \tau \gamma u_k^n v_k^n q^{n-\frac{1}{2}} \\
		&{}+ \tau \gamma u_k^n v_k^n q^{n+\frac{1}{2}} =  \tau \alpha_1 \hat{u}^n 
		+ 2 \tau \gamma u_k^n v_k^n \left [ (q_k^{n-\frac{1}{2}} - p_k^{n-\frac{1}{2}}) +
		(q_k^{n+\frac{1}{2}} - p_k^{n+\frac{1}{2}}) \right ]
	\end{aligned}
	\label{eq:SQP_adjointp_semi}
\end{equation}
and
\begin{equation}
	\begin{aligned}
		& \tau \gamma \left [ (q_k^{n-\frac{1}{2}} - p_k^{n-\frac{1}{2}}) +
		(q_k^{n+\frac{1}{2}} - p_k^{n+\frac{1}{2}}) \right ] u_k^n u^n 
		+ \tau \alpha_2 v^n - 
		\frac{\tau \gamma}{2}  ( u_k^n )^2 p^{n-\frac{1}{2}}
		\\ 
		&{} -\frac{\tau \gamma }{2} ( u_k^n )^2 p^{n+\frac{1}{2}}  
		+ q^{n-\frac{1}{2}}
		- \frac{\tau D_v}{2}\nabla^2 q^{n-\frac{1}{2}} 
		+ \frac{\tau \gamma}{2}  \left(u_k^n \right)^2  q^{n-\frac{1}{2}} 
		- q^{n+\frac{1}{2}} \\ 
		&{}- \frac{\tau D_v}{2}\nabla^2 q^{n+\frac{1}{2}} 
		  +\frac{\tau \gamma }{2} \left(u_k^n \right)^2  q^{n+\frac{1}{2}} = 		
		\tau \alpha_2 \hat{v}^n
		+ \tau \gamma \left(u_k^n \right)^2 \left [ (q_k^{n-\frac{1}{2}} - p_k^{n-\frac{1}{2}}) \right. \\
            &{}+\left.(q_k^{n+\frac{1}{2}} - p_k^{n+\frac{1}{2}}) \right ].
	\end{aligned}
	\label{eq:SQP_adjointq_semi}
\end{equation}
At the final half-time-step, we have
\begin{equation}
	\begin{aligned}
		&\frac{\tau \alpha_1}{2} u^{N_t} + \tau \gamma ( q_k^{N_t-\frac{1}{2}} - p_k^{N_t-\frac{1}{2}}) v_k^{N_t} u^{N_t}
  + \tau \gamma  (q_k^{N_t-\frac{1}{2}} - p_k^{N_t-\frac{1}{2}}) u_k^{N_t} v^{N_t} 
		+ p^{N_t-\frac{1}{2}} \\
        &{} - \frac{\tau D_u}{2}\nabla^2 p^{N_t-\frac{1}{2}} + \frac{\tau \gamma }{2} \left(1- 2 u_k^{N_t} v_k^{N_t}\right)  p^{N_t-\frac{1}{2}}  + \tau \gamma u_k^{N_t} v_k^{N_t} q^{N_t-\frac{1}{2}}
		= \frac{\tau \alpha_1}{2} \hat{u}^{N_t} \\
		&{}+ 2 \tau \gamma u_k^{N_t} v_k^{N_t} (q_k^{N_t-\frac{1}{2}} - p_k^{N_t-\frac{1}{2}})
	\end{aligned}
	\label{eq:SQP_adjointp_semi_pNt}
\end{equation}
and
\begin{equation}
	\begin{aligned}
		& \tau \gamma (q_k^{N_t-\frac{1}{2}} - p_k^{N_t-\frac{1}{2}})u_k^{N_t} u^{N_t} 
		+ \frac{\tau \alpha_2}{2} v^{N_t} - 
		\frac{\tau \gamma}{2}  ( u_k^{N_t} )^2 p^{N_t-\frac{1}{2}}
		+ q^{N_t-\frac{1}{2}} \\ 
		&{}- \frac{\tau D_v}{2}\nabla^2 q^{N_t-\frac{1}{2}} 
		+ \frac{\tau \gamma}{2}  (u_k^{N_t} )^2  q^{N_t-\frac{1}{2}} 
		 = 	\frac{\tau \alpha_2}{2} \hat{v}^{N_t}
		+ \tau \gamma (u_k^{N_t} )^2 (q_k^{N_t-\frac{1}{2}} - p_k^{N_t-\frac{1}{2}}).
	\end{aligned}
	\label{eq:SQP_adjointq_semi_qNt}
\end{equation}

\subsection{Finite element discretization}
\label{sec:otd_finite_elements}
For the spatial discretization of \eqref{eq:SQP_stateu_semi}--\eqref{eq:SQP_adjointq_semi_qNt} we apply a finite element method 
and use P$1$ elements and $N_x$ nodes in the discretization of the domain $\Omega$. 
Let $\mathbf{M}$ denote the mass matrix and $\mathbf{K}$ the stiffness matrix, with entries
$$\left[\mathbf{M}\right]_{r,s} = \int_\Omega \varphi_r \varphi_s \, dx \; \; \text{ and } \; \; \left[\mathbf{K}\right]_{r,s} = \int_\Omega \nabla \varphi_r \cdot \nabla \varphi_s \, dx, $$
where $\varphi_r$, $\varphi_s$ denote the finite element basis functions. Both these matrices are symmetric, $\mathbf{M}$ is positive definite, and $\mathbf{K}$ is positive semi-definite (positive definite when considering Dirichlet boundary conditions). For the terms arising from linearization, we define the pseudo-mass matrices $\mathbf{M}_{u^i v^i}$,  $\mathbf{M}_{u^i(q^j - p^j)}$, $\mathbf{M}_{(u^i)^2}$, and $ \mathbf{M}_{v^i(q^j - p^j)}$ as 
\begin{equation*}
	\begin{aligned}
		\left[\mathbf{M}_{u^i v^i}\right]_{r,s} := \int_\Omega u_k^i v_k^i \varphi_r \varphi_s \, dx,  
		\hspace{5mm}
		&\left[ \mathbf{M}_{u^i(q^j - p^j)}\right] _{r,s} := \int_\Omega u_k^i(q_k^j - p_k^j) \varphi_r \varphi_s \, dx, \\
		\left[ \mathbf{M}_{(u^i)^2}\right] _{r,s} := \int_\Omega(u_k^i)^2 \varphi_r \varphi_s \, dx, 
		\hspace{5mm}
		&\left[ \mathbf{M}_{v^i(q^j - p^j)}\right] _{r,s} := \int_\Omega v_k^i(q_k^j - p_k^j) \varphi_r \varphi_s \, dx,
	\end{aligned}
\end{equation*}
which are also symmetric. 
The superscripts $i$ and $j$ of the variables involved denote the $i$-th and $j$-th time points ($t^i=i \tau$ and $t^j=j \tau$). 
For the right-hand sides of the discretized equations, we define the vectors $\bm{d}^i$, $\bm{c}^{i,j}$ and $\bm{h}^{i,j}$, given by 
\begin{equation*}
	\begin{aligned}
	\left[  \bm{d}^i \right]_{r} &:= \int_\Omega 2 \gamma (u_k^i)^2 v_k^i \varphi_r \, dx, \\
	\left[  \bm{c}^{i,j} \right]_{r} &:=  \int_\Omega \left[ \alpha_1 \hat{u}^i + 4\gamma u_k^i v_k^i  (q_k^j - p_k^j ) \right] \varphi_r \, dx, \\
	\left[  \bm{h}^{i,j} \right]_{r} &:= \int_\Omega \left[ \alpha_2 \hat{v}^i + 2\gamma ( u_k^i)^2 (q_k^j - p_k^j )\right] \varphi_r \, dx,
	\end{aligned}
\end{equation*}
where  $\hat{u}^i(\textbf{x})=\hat{u}(t^i,\textbf{x})$ and $\hat{v}^i(\textbf{x})=\hat{v}(t^i,\textbf{x})$. 
We also define $\bm{\hat{u}}^i$ with $[\bm{\hat{u}}^i]_r = \int_\Omega \hat{u}^i \psi_r \, dx$, similarly for $\bm{\hat{v}}^i$.
Incorporating the above matrices and vectors yields the fully discretized state equations
\begin{equation}
	\begin{aligned}
		&\left[-\mathbf{M} + \frac{\tau D_u}{2} \mathbf{K} + \frac{\tau \gamma}{2} \mathbf{M} - \tau \gamma \mathbf{M}_{u^n v^n} \right] \bm{u}^n 
		+ \left[ \mathbf{M} + \frac{\tau D_u}{2} \mathbf{K} + \frac{\tau \gamma}{2} \mathbf{M} \right.\\
            &{}-\left. \tau \gamma \mathbf{M}_{u^{n+1} v^{n+1}} \right] \bm{u}^{n+1} 
		- \frac{\tau \gamma}{2} \mathbf{M}_{(u^n)^2 } \bm{v}^n 
		- \frac{\tau \gamma}{2} \mathbf{M}_{(u^{n+1})^2 } \bm{v}^{n+1} \\
		&{}- \frac{\tau \gamma}{2} \mathbf{M} \bm{a}^{n+\frac{1}{2}} = - \frac{\tau}{2} \left( \bm{d}^n + \bm{d}^{n+1} \right) 
	\end{aligned}
	\label{eq:otd_stateu_discr}
\end{equation}
and 
\begin{equation}
	\begin{aligned}
		&\tau \gamma \mathbf{M}_{u^n v^n} \bm{u}^n + \tau \gamma \mathbf{M}_{u^{n+1} v^{n+1}} \bm{u}^{n+1} +
		\left[-\mathbf{M} + \frac{\tau D_v}{2} \mathbf{K} + \frac{\tau \gamma}{2} \mathbf{M}_{(u^n)^2 } \right] \bm{v}^n  \\
		& {}+ \left[ \mathbf{M} + \frac{\tau D_v}{2} \mathbf{K} + \frac{\tau \gamma}{2} \mathbf{M}_{(u^{n+1})^2}\right] \bm{v}^{n+1} 
		- \tau \gamma \mathbf{M} \bm{b}^{n+\frac{1}{2}} = \frac{\tau}{2} \left( \bm{d}^n + \bm{d}^{n+1} \right),
	\end{aligned}
	\label{eq:otd_statev_discr}
\end{equation}
with initial conditions $\mathbf{M} \bm{u}^0 = \bm{u}_0$, $\mathbf{M} \bm{v}^0 = \bm{v}_0$, where $[\bm{u}_0]_r = \int_\Omega u_0 \psi_r \, dx$, similarly for $\bm{v}_0$, using the initial conditions given in \eqref{eq:schnakenberg_BC}.
For the first half step, the discretized adjoint equations are
\begin{equation*}
	\begin{aligned}
		&\frac{\tau \alpha_1}{2} \mathbf{M} \bm{u}^0 + \tau \gamma \mathbf{M}_{v^0(q^{\frac{1}{2}} - p^{\frac{1}{2}})}  \bm{u}^0 
        + \tau \gamma 
		\mathbf{M}_{u^0(q^{\frac{1}{2}} - p^{\frac{1}{2}})} \bm{v}^0 
		+ \mathbf{M} \bm{p}^0 \\
  & {}+ \left[- \mathbf{M} + \frac{\tau D_u}{2} \mathbf{K} + \frac{\tau \gamma}{2} \mathbf{M} - \tau \gamma \mathbf{M}_{u^0v^0} \right] \bm{p}^{\frac{1}{2}}    
		+ \tau \gamma \mathbf{M}_{u^0 v^0} \bm{q}^{\frac{1}{2}} = \frac{\tau}{2}\bm{c}^{0,\frac{1}{2}} \\
	\end{aligned}
\end{equation*}
and
\begin{equation*}
\begin{aligned} 
	 &\tau \gamma  \mathbf{M}_{u^0(q^{\frac{1}{2}} - p^{\frac{1}{2}})} \bm{u}^0 + \frac{\tau \alpha_2}{2} \mathbf{M} \bm{v}^0 
	 -
	\frac{\tau \gamma}{2} \mathbf{M}_{(u^0)^2} \bm{p}^{\frac{1}{2}} 
	+ \mathbf{M} \bm{q}^0 \\
 & {} + \left[-\mathbf{M} + \frac{\tau D_v}{2} \mathbf{K} + \frac{\tau \gamma}{2} \mathbf{M}_{(u^0)^2} \right] \bm{q}^{\frac{1}{2}} 
	= \frac{\tau}{2} \bm{h}^{0,\frac{1}{2}}.
\end{aligned}
\end{equation*}
For $n=1,...,N_t-1$, we obtain
\begin{equation}
	\begin{aligned}
		&\tau \alpha_1 \mathbf{M} \bm{u}^n 
            + \tau \gamma \left[ \mathbf{M}_{v^n(q^{n-\frac{1}{2}} - p^{n-\frac{1}{2}})} 
		+ \mathbf{M}_{v^n(q^{n+\frac{1}{2}} - p^{n+\frac{1}{2}})} \right] \bm{u}^n \\
            &{}+ \tau \gamma \left[ \mathbf{M}_{u^n(q^{n-\frac{1}{2}} - p^{n-\frac{1}{2}})} + \mathbf{M}_{u^n(q^{n+\frac{1}{2}} - p^{n+\frac{1}{2}})} \right] \bm{v}^n 
            + \left[ \mathbf{M} + \frac{\tau D_u}{2} \mathbf{K} + \frac{\tau \gamma}{2} \mathbf{M} \right.\\
            &\left. {}- \tau \gamma \mathbf{M}_{u^nv^n} \right] \bm{p}^{n-\frac{1}{2}}
		+ \left[- \mathbf{M} + \frac{\tau D_u}{2} \mathbf{K} +  \frac{\tau \gamma}{2} \mathbf{M} - \tau \gamma \mathbf{M}_{u^nv^n} \right] \bm{p}^{n+\frac{1}{2}} \\
		&{}+ \tau \gamma \mathbf{M}_{u^n v^n} \bm{q}^{n-\frac{1}{2}} 
		+ \tau \gamma \mathbf{M}_{u^n v^n} \bm{q}^{n+\frac{1}{2}}  = \frac{\tau}{2} \left( \bm{c}^{n, n-\frac{1}{2}} + \bm{c}^{n, n+\frac{1}{2}} \right)
	\end{aligned}
	\label{eq:otd_adjp_discr}
\end{equation}
and
\begin{equation}
\begin{aligned} 
	&\tau \gamma \left[ \mathbf{M}_{u^n(q^{n-\frac{1}{2}} - p^{n-\frac{1}{2}})} + \mathbf{M}_{u^n(q^{n+\frac{1}{2}} - p^{n+\frac{1}{2}})} \right] \bm{u}^n 
	+\tau \alpha_2 \mathbf{M} \bm{v}^n 
	- \frac{\tau \gamma}{2} \mathbf{M}_{(u^n)^2} \bm{p}^{n-\frac{1}{2}} \\
	&-
	\frac{\tau \gamma}{2} \mathbf{M}_{(u^n)^2} \bm{p}^{n+\frac{1}{2}}    
	+ \left[ \mathbf{M}  + \frac{\tau D_v}{2} \mathbf{K} + \frac{\tau \gamma}{2} \mathbf{M}_{(u^n)^2} \right] \bm{q}^{n-\frac{1}{2}} \\
	&+ \left[-\mathbf{M} + \frac{\tau D_v}{2} \mathbf{K} + \frac{\tau \gamma}{2} \mathbf{M}_{(u^n)^2} \right] \bm{q}^{n+\frac{1}{2}} 
	= \frac{\tau}{2} \left( \bm{h}^{n, n-\frac{1}{2}} + \bm{h}^{n, n+\frac{1}{2}} \right),
\end{aligned}
	\label{eq:otd_adjq_discr}
\end{equation}
and at the final half step, we have
\begin{equation*}
	\begin{aligned}
		&\tau \alpha_1 \mathbf{M} \bm{u}^{N_t} + \tau \gamma  \mathbf{M}_{v^{N_t}(q^{N_t-\frac{1}{2}} - p^{N_t-\frac{1}{2}})} 
		\bm{u}^{N_t} + 
		\tau \gamma  \mathbf{M}_{u^{N_t}(q^{N_t-\frac{1}{2}} - p^{N_t-\frac{1}{2}})}  \bm{v}^{N_t}  
		+ \left[ \mathbf{M} \right.\\ 
            &\left. {}+ \frac{\tau D_u}{2} \mathbf{K} + \frac{\tau \gamma}{2} \mathbf{M} - \tau \gamma \mathbf{M}_{u^{N_t}v^{N_t}} \right] \bm{p}^{N_t-\frac{1}{2}} 
		+ \tau \gamma \mathbf{M}_{u^{N_t} v^{N_t}} \bm{q}^{N_t-\frac{1}{2}}   
	       = \frac{\tau}{2}\bm{c}^{N_t,N_t-\frac{1}{2}} \\
	\end{aligned}
\end{equation*}
and 
\begin{equation*}
\begin{aligned} 
	&\tau \gamma \mathbf{M}_{u^{N_t}(q^{N_t-\frac{1}{2}} - p^{N_t-\frac{1}{2}})} \bm{u}^{N_t} 
	+\tau \alpha_2 \mathbf{M} \bm{v}^{N_t} 
	- \frac{\tau \gamma}{2} \mathbf{M}_{(u^{N_t})^2} \bm{p}^{N_t-\frac{1}{2}} \\
	& {}+ \left[ \mathbf{M}  + \frac{\tau D_v}{2} \mathbf{K} + \frac{\tau \gamma}{2} \mathbf{M}_{(u^{N_t})^2} \right] \bm{q}^{N_t-\frac{1}{2}} = \frac{\tau}{2} \bm{h}^{N_t,N_t-\frac{1}{2}}.
\end{aligned}
\end{equation*}

Finally, according to the gradient equations \eqref{eq:OTD_gradient}, the controls depend on the adjoint variables which we evaluate at half time-steps:
\begin{equation}
	\begin{aligned}
		\beta_1 \bm{a}^{n+\frac{1}{2}} - \gamma \bm{p}^{n+\frac{1}{2}}  = \bm{0}, \hspace{5mm} 
		\beta_2 \bm{b}^{n+\frac{1}{2}} - \gamma \bm{q}^{n+\frac{1}{2}}  = \bm{0}, \quad \text{ for } \; n=0,1,...,N_t-1.
	\end{aligned}
	\label{eq:OTD_gradient_discr}
\end{equation}
We have thus arrived at a fully discretized system following the Optimize-then-Discretize approach. To solve the system, we will form a system matrix for the all-at-once problem. First, however, 
we show that the St\"{o}rmer--Verlet method leads to the same discrete optimality conditions when considering the Discretize-then-Optimize approach.

\section{Discretize-then-Optimize approach}
\label{sec:dto}
In this section, we approximate the Lagrangian for the SQP problem by introducing a time-stepping strategy treating each term individually, to achieve a match with the time discretization in Section~\ref{sec:otd}.  The difficulty arises from the linearized terms that couple the state and adjoint variables, thus combining integer and half time-steps within the same term.

We note that the  St\"{o}rmer--Verlet scheme \eqref{eq:stormer_verlet2}  is similar to a trapezoidal rule approximation, equivalent to the midpoint rule for linear equations. We therefore use  these quadrature rules,  with the notation $u^n \approx u(t^n)$ and $u^{n+\frac{1}{2}} \approx u(t^{n+\frac{1}{2}})$, where $t^{n+\frac{1}{2}}:= t^n + \frac{\tau}{2}$. The latter term can in turn be approximated to second order by $\frac{1}{2} \left(u^n + u^{n+1}\right)$. 

We first examine the SQP cost functional $\mathcal{\tilde{J}}$ defined in \eqref{eq:SQP_Jtilde_final}.
Starting with the integrals involving $\left( u-\hat{u}\right)^2$ and $\left( v-\hat{v}\right)^2$, we can apply the trapezoidal quadrature rule to obtain 
\begin{equation*}
	\frac{\alpha_1}{2}  \int_Q \left( u-\hat{u}\right)^2 \, d Q  \approx \frac{\tau \alpha_1}{4}  \int_\Omega \Big[(u^0 - \hat{u}^0)^2 + 
	2\sum_{i=1}^{N_t-1}(u^i - \hat{u}^i)^2 +
	(u^{N_t} - \hat{u}^{N_t})^2  \Big] \, dx,
\end{equation*}
and similarly for $\left( v-\hat{v}\right)^2$. The control and adjoint variables need to be discretized at half steps, i.e., at $\frac{1}{2},...,N_t-\frac{1}{2}$. For the control variables, we apply the midpoint rule to the integrals involving $a^2$ and $b^2$:
\begin{equation*}
	\frac{\beta_1}{2} \int_Q  a^2 \, d Q  \approx \frac{\tau \beta_1}{2} \int_\Omega  
	\sum_{i=0}^{N_t-1} (a^{i+\frac{1}{2}})^2  \, dx, 
\end{equation*}
and similarly for $b^2$. 
The rest of the terms are integrals of polynomials involving adjoint variables coupled with state variables. Here, we also apply the midpoint rule, and to achieve an analogous time-stepping to the OTD approach, we apply a further second-order approximation in time. We illustrate the application of these steps on the following term, with other calculations following similarly:
 \begin{equation*}
 	\begin{aligned}
 		\int_Q &\gamma \left( q_k - p_k\right)  v_k \left( u-u_k\right)^2 \, d Q  \\
 		&\approx  \int_\Omega \tau \gamma \sum_{i=0}^{N_t-1} (q_k^{i+\frac{1}{2}} - p_k^{i+\frac{1}{2}})v_k^{i+\frac{1}{2}} (u^{i+\frac{1}{2}}-u_k^{i+\frac{1}{2}} )^2 \, dx \\
 		&\approx 
 		\int_\Omega \frac{\tau\gamma}{2} \left ( 
 		(q_k^{\frac{1}{2}} - p_k^{\frac{1}{2}}) 
 		v_k^0 ( u^0-u_k^0 ) ^2 + 
 		\sum_{i=1}^{N_t-1} 
 		\left[ (q_k^{i-\frac{1}{2}} - p_k^{i-\frac{1}{2}}) \right.\right. \\
 		&\hspace{3mm} {}+ \left.\left. (q_k^{i+\frac{1}{2}} - p_k^{i+\frac{1}{2}}) \right]
		 v_k^i (u^i-u_k^i )^2 
		 + (q_k^{N_t-\frac{1}{2}} - p_k^{N_t-\frac{1}{2}}) 
 		v_k^{N_t} ( u^{N_t}-u_k^{N_t} ) ^2 
 		\right) \, dx.
	\end{aligned}
\end{equation*}

Next, we approximate the constraints using a similar two-step scheme as the one above. 
For instance, applying the midpoint method to the first state equation in \eqref{eq:SQP_state} yields 
\begin{flalign*}
		u^{n+1} - u^n - \tau D_u \nabla^2 u^{n+\frac{1}{2}} + \tau \gamma \big( 1-2u_k^{n+\frac{1}{2}}  v_k^{n+\frac{1}{2}}\big) u^{n+\frac{1}{2}}  - \tau \gamma (u_k^{n+\frac{1}{2}} )^2 v^{n+\frac{1}{2}}  \\ 
		-\tau \gamma a^{n+\frac{1}{2}} + 2 \tau \gamma (u_k^{n+\frac{1}{2}} )^2 v_k^{n+\frac{1}{2}} &\approx 0, 
\end{flalign*}
for $n=0,1,...,N_t-1$. Then, approximating the terms that involve state variables at half-time-steps using the trapezoidal rule, we obtain
\begin{equation*}
	\begin{aligned}
		&u^{n+1} - u^n 
		- \frac{\tau D_u}{2} \nabla^2 \left(u^n + u^{n+1}\right)  
		+ \frac{\tau \gamma}{2} \left[ (1-2u_k^n  v_k^n)u^n \right. \\
            &+\left. (1-2u_k^{n+1}  v_k^{n+1})u^{n+1}\right] 
		-\frac{\tau \gamma}{2} \left[ \left(u_k^n \right)^2 v^n + \left(u_k^{n+1}  \right)^2 v^{n+1} \right] \\
		&-\tau \gamma a^{n+\frac{1}{2}}
		 + \tau \gamma \left[ \left(u_k^n \right)^2 v_k^n + \left(u_k^{n+1} \right)^2 v_k^{n+1}\right] \approx 0.
	\end{aligned}
\end{equation*}
This leads to integer time-steps for the state variables and half time-steps for the control variables. For each $n$ we multiply the discretized equation by the multiplier $p^{n+\frac{1}{2}}$, 
and add the integral of this term to the cost functional. We repeat the same process with the second equation in~\eqref{eq:SQP_state} with the corresponding multiplier $q^{n+\frac{1}{2}}$. 
Then we discretize the Lagrangian using the finite element method and obtain the fully discrete problem with the mass, stiffness, and pseudo-mass matrices as  in Section~\ref{sec:otd_finite_elements}.

To obtain the optimality conditions from the discrete Lagrangian, we compute its partial derivatives. Differentiating with respect to the adjoint variables $\bm{p}^{\frac{1}{2}}, ..., \bm{p}^{N_t-\frac{1}{2}}$ and $\bm{q}^{\frac{1}{2}}, ..., \bm{q}^{N_t-\frac{1}{2}}$ gives the same system of discretized state equations as 	\eqref{eq:otd_stateu_discr} and \eqref{eq:otd_statev_discr}, respectively.
Differentiation with respect to $\bm{u}^1, ..., \bm{u}^{N_t}$ and $\bm{v}^1, ..., \bm{v}^{N_t}$ yields  adjoint equations \eqref{eq:otd_adjp_discr}--\eqref{eq:otd_adjq_discr}, and we also recover the gradient equations as in \eqref{eq:OTD_gradient_discr}. 

\section{All-at-once system}
\label{sec:all_at_once}
We have constructed a methodology by which the linearized systems arising from the OTD and DTO approaches coincide. To solve for the state, control, and adjoint variables, we can now construct the all-at-once system matrix.  
For a more convenient notation, we introduce the matrices
\begin{equation*}
	\begin{aligned}
		\mathbf{L}_{(1)}^i := \frac{D_u}{2} \mathbf{K} + \frac{\gamma}{2} \mathbf{M} - \gamma \mathbf{M}_{u^iv^i}, \hspace{2mm}
		\mathbf{L}_{(2)}^i := \frac{D_v}{2} \mathbf{K} + \frac{\gamma}{2}             \mathbf{M}_{(u^i)^2},  \hspace{2mm}
            \text{ for } \; i=0,...,N_t,
	\end{aligned}
\end{equation*}
and
\begin{equation*}
	\begin{aligned}
		\mathbf{A}_{(1)}^i := \begin{cases}
			 \gamma \mathbf{M}_{v^0(q^\frac{1}{2} - p^\frac{1}{2})}, & i=0, \\
			  \gamma \left (\mathbf{M}_{v^i(q^{i-\frac{1}{2}} - p^{i-\frac{1}{2}})} + \mathbf{M}_{v^i(q^{i+\frac{1}{2}} - p^{i+\frac{1}{2}})} \right), & i=1,...,N_t-1, \\
			  \gamma \mathbf{M}_{v^{N_t}(q^{N_t-\frac{1}{2}} - p^{N_t-\frac{1}{2}})}, & i=N_t,
		\end{cases} 
		\\
		\mathbf{A}_{(12)}^i := \begin{cases}
			\gamma \mathbf{M}_{u^0(q^\frac{1}{2} - p^\frac{1}{2})}, & i=0, \\
			\gamma \left(\mathbf{M}_{u^i(q^{i-\frac{1}{2}} - p^{i-\frac{1}{2}})} + \mathbf{M}_{u^i(q^{i+\frac{1}{2}} - p^{i+\frac{1}{2}})} \right), & i=1,...,N_t-1, \\
		\gamma \mathbf{M}_{u^{N_t}(q^{N_t-\frac{1}{2}} - p^{N_t-\frac{1}{2}})}, & i=N_t.
		\end{cases}
	\end{aligned}
\end{equation*}
To devise an efficient iterative solver, we transform the system matrix into a simpler and more convenient form. 
Firstly, using the initial conditions for the state equations, we can separate the adjoint equations at the first time-step and obtain 
\begin{equation}
	\begin{aligned}
		\mathbf{M} \bm{p}^0 &= \frac{\tau \alpha_1}{2} \left(\bm{\hat{u}}^0 - \bm{u}_0  \right) - \left[ -\mathbf{M} + \tau \mathbf{L}_{(1)}^0 \right] \bm{p}^\frac{1}{2} - \tau \gamma \mathbf{M}_{u^0v^0} \bm{q}^\frac{1}{2}, \\
		\mathbf{M} \bm{q}^0 &= \frac{\tau \alpha_2}{2} \left(\bm{\hat{v}}^0 - \bm{v}_0  \right) + \frac{\tau \gamma }{2} \mathbf{M}_{\left( u^0\right)^2}\bm{p}^\frac{1}{2}  - \left[ -\mathbf{M} + \tau \mathbf{L}_{(2)}^0 \right]\bm{q}^\frac{1}{2},
	\end{aligned}
        \label{eq:aat_pq0}
\end{equation}
for $\bm{p}^0$ and $\bm{q}^0$. 
In block notation, the remaining all-at-once system can be written as 
\begin{equation}\label{allatones_system}
	\begin{bmatrix}
		\mathbf{A}_1 & \mathbf{0} & \mathbf{B}_1^\top \\
		\mathbf{0} & \mathbf{A}_2 & \mathbf{B}_2^\top \\
		\mathbf{B}_1 & \mathbf{B}_2 & \mathbf{0}
	\end{bmatrix}
	\begin{bmatrix}
		\underline{\bm{u}} \\ \underline{\bm{a}} \\ \underline{\bm{p}}
	\end{bmatrix} = 
	\begin{bmatrix}
		\underline{\bm{c}} \\ \underline{\bm{0}} \\ \underline{\bm{d}}
	\end{bmatrix},
\end{equation}
were the vectors for the fully discretized state, adjoint and control variables are defined as 
\begin{equation*}
	\begin{aligned}
			\underline{\bm{u}} &:= \left[(\bm{u}^1)^{\top}, ..., (\bm{u}^{N_t})^{\top}, (\bm{v}^1)^{\top}, ..., (\bm{v}^{N_t})^{\top}\right]^{\top}, \\
			\underline{\bm{p}} &:= \left[(\bm{p}^{\frac{1}{2}})^{\top}, ..., (\bm{p}^{N_t-\frac{1}{2}})^{\top}, (\bm{q}^{\frac{1}{2}})^{\top}, ..., (\bm{q}^{N_t-\frac{1}{2}})^{\top}\right]^{\top},\\
			\underline{\bm{a}} &:= \left[(\bm{a}^{\frac{1}{2}})^{\top}, ..., (\bm{a}^{N_t-\frac{1}{2}})^{\top}, (\bm{b}^{\frac{1}{2}})^{\top}, ..., (\bm{b}^{N_t-\frac{1}{2}})^{\top}\right]^{\top},
		\end{aligned}
\end{equation*}
and the right-hand side vectors are defined by the right-hand side terms derived in Section~\ref{sec:otd_finite_elements}, with $\underline{\bm{c}}$ and $\underline{\bm{d}}$ corresponding to the adjoint and state equations, respectively. The block-matrices in \eqref{allatones_system} are
\begin{equation*}
	\begin{aligned}
			\mathbf{A}_1 :=& \begin{bmatrix}
			\mathbf{A}_{11} & \mathbf{A}_{12} \\
			\mathbf{A}_{21} & \mathbf{A}_{22}
		\end{bmatrix}, 
		\hspace{5mm}
		\mathbf{A}_2:=\text{blkdiag}\big( \left\lbrace \tau \beta_1 \mathbf{M}\right\rbrace_{i=1}^{N_t}  , \left\lbrace  \tau \beta_2 \mathbf{M}\right\rbrace_{i=1}^{N_t}   \big) , \\
		\mathbf{B}_1 :=& \begin{bmatrix}
			\mathbf{B}_{11} & \mathbf{B}_{12} \\
			\mathbf{B}_{21} & \mathbf{B}_{22}
		\end{bmatrix},
		\hspace{5mm}
		\mathbf{B}_2:=\text{blkdiag}\big( \left\lbrace -\tau \gamma \mathbf{M}\right\rbrace_{i=1}^{2N_t}\big).
	\end{aligned} 
\end{equation*}
The sub-blocks of $\mathbf{A}_1$ are
\begin{equation*}
	\begin{aligned}
	\mathbf{A}_{11} &:= \text{blkdiag} \left( \left\lbrace \tau \alpha_1 \mathbf{M} + \tau \mathbf{A}_{(1)}^i \right\rbrace _{i=1}^{N_t-1}, \frac{\tau \alpha_1}{2} \mathbf{M} + \tau \mathbf{A}_{(1)}^{N_t}   \right), \\
	\mathbf{A}_{12} \equiv \mathbf{A}_{21} &:= \text{blkdiag} \left( \left\lbrace  \tau \mathbf{A}_{(12)}^i \right\rbrace _{i=1}^{N_t}  \right), \\
	\mathbf{A}_{22} &:= \text{blkdiag} \left( \left\lbrace \tau \alpha_2 \mathbf{M} \right\rbrace _{i=1}^{N_t-1}, \frac{\tau \alpha_2}{2} \mathbf{M}   \right), 
	\end{aligned}
\end{equation*}
while the (lower-triangular) sub-blocks of $\mathbf{B}_1$ are
\begin{equation*}
	\begin{aligned}
		\mathbf{B}_{11} &:= \text{blkdiag} \left( \left\lbrace \mathbf{M} + \tau \mathbf{L}_{(1)}^i \right\rbrace _{i=1}^{N_t} \right) + \text{blksubdiag} \left(\left\lbrace -\mathbf{M} + \tau \mathbf{L}_{(1)}^i\right\rbrace_{i=1}^{N_t-1}  \right),  \\
		\mathbf{B}_{12} &:= \text{blkdiag} \left( \left\lbrace -\frac{\tau\gamma}{2} \mathbf{M}_{\left( u^i\right)^2} \right\rbrace _{i=1}^{N_t} \right) + \text{blksubdiag} \left(\left\lbrace  -\frac{\tau\gamma}{2} \mathbf{M}_{\left( u^i\right)^2} \right\rbrace_{i=1}^{N_t-1}  \right),  \\
		\mathbf{B}_{21} &:= \text{blkdiag} \left( \left\lbrace \tau\gamma \mathbf{M}_{u^iv^i} \right\rbrace _{i=1}^{N_t} \right) + \text{blksubdiag} \left(\left\lbrace  \tau\gamma \mathbf{M}_{u^iv^i} \right\rbrace_{i=1}^{N_t-1}  \right),  \\
		\mathbf{B}_{22} &:= \text{blkdiag} \left( \left\lbrace \mathbf{M} + \tau \mathbf{L}_{(2)}^i \right\rbrace _{i=1}^{N_t} \right) + \text{blksubdiag} \left(\left\lbrace -\mathbf{M} + \tau \mathbf{L}_{(2)}^i \right\rbrace_{i=1}^{N_t-1}  \right).
	\end{aligned}
\end{equation*}

We can now eliminate the discretized gradient equations and reorder the equations and variables.
With two state and two adjoint PDEs, we can think of the system matrix as of a matrix composed of four blocks. For computational convenience, we wish for the $(1,1)$-block of the matrix to be positive semi-definite. This can be implemented by absorbing the negative sign into part of the vector of unknowns
and by negating the $(2,2)$-block of the matrix and part of the right-hand side.
With these modifications, we obtain the system
\begin{equation}
    \underbrace{\begin{bmatrix}
			\mathbf{A} &  \mathbf{B}^\top \\
			\mathbf{B} & -\mathbf{C}
	\end{bmatrix}}_{\bm{\mathcal{A}}}
    \underbrace{\begin{bmatrix}
		-\underline{\bm{p}} \\ \underline{\bm{u}}
	\end{bmatrix}}_{\underline{\bm{w}}} := 
    \begin{bmatrix}
		\mathbf{A}_3 &  \mathbf{B}_1 \\
		\mathbf{B}_1^\top & -\mathbf{A}_1
	\end{bmatrix}
	\begin{bmatrix}
		-\underline{\bm{p}} \\ \underline{\bm{u}}
	\end{bmatrix} = 
    \underbrace{\begin{bmatrix}
		\underline{\bm{d}} \\ -\underline{\bm{c}}
	\end{bmatrix}}_{\underline{\bm{b}}},
	\label{eq:allatonce_sys}
\end{equation}
where the block $\mathbf{A}_3:= \mathbf{B}_2 \mathbf{A}_2^{-1} \mathbf{B}_2^\top$ arises by incorporating the gradient equations:
\begin{equation*}
	\mathbf{A}_3:=\text{blkdiag}\left(  \left\lbrace \frac{\tau \gamma ^2}{\beta_1} \mathbf{M}\right\rbrace_{i=1}^{N_t}  , \left\lbrace  \frac{\tau \gamma ^2}{\beta_2}\mathbf{M}\right\rbrace_{i=1}^{N_t} \right).
\end{equation*}
We highlight that the system matrix $\bm{\mathcal{A}}$ is symmetric indefinite and the main challenge  is to find a solution to the large-scale (generalized) saddle-point problem~\eqref{eq:allatonce_sys}. 
After solving the system above
and simplifying the terms, we solve \eqref{eq:aat_pq0} for $\bm{p}^0$ and $\bm{q}^0$
using Chebyshev semi-iteration~\cite{golubvarga1961,wathenrees2009} with Jacobi splitting. 

\section{Preconditioning} \label{sec:precond}

Solving large (generalized) saddle-point problems of the form \eqref{eq:allatonce_sys} is a topic of wide interest in the numerical analysis community. Since finding an ``exact'' solution to this problem by using a direct solver is computationally too expensive, we make use of iterative solvers known to be suitable for this type of problem structure, specifically  Krylov subspace solvers. We need to ensure that we find the solutions of the systems robustly and efficiently, that the iterative solver will converge fast and we take advantage of the problem structures within the algorithms used for their solution.

To address the first issue, we note that the convergence rate of iterative methods often depends on the distribution of eigenvalues of matrix $\bm{\mathcal{A}}$. As $\bm{\mathcal{A}}$ is considered to be indefinite, and its eigenvalues can be distributed across many orders of magnitude, the convergence of an iterative method can be slow.  By constructing a preconditioner $\bm{\bm{\mathcal{P}}}$, we can transform the linear system \eqref{eq:allatonce_sys} into the left-preconditioned system
$\bm{\mathcal{P}}^{-1} \bm{\mathcal{A}}\underline{\bm{w}} = \bm{\mathcal{P}}^{-1}\underline{\bm{b}}$ or the right-preconditioned system
$(\bm{\mathcal{A}}\bm{\mathcal{P}}^{-1}) (\bm{\mathcal{P}}\underline{\bm{w}}) = \underline{\bm{b}}$,  where $\bm{\mathcal{P}}^{-1} \bm{\mathcal{A}}$ (equivalently $\bm{\mathcal{A}}\bm{\mathcal{P}}^{-1}$) has improved properties that ensure faster convergence of iterative methods. 
An ``ideal'' preconditioner for (generalized) saddle-point systems~\eqref{eq:allatonce_sys}, given certain properties of matrices $\mathbf{A}$, $\mathbf{B}$ and $\mathbf{C}$, is of the form
\begin{equation}
	\bm{\mathcal{P}}= \begin{bmatrix}
		\mathbf{A} & \mathbf{0} \\
		\mathbf{0} & \mathbf{S}
	\end{bmatrix},
	\label{eq:ideal_prec}
\end{equation}
where $\mathbf{S}$ denotes the (negative) Schur complement $\mathbf{S} := \mathbf{C} + \mathbf{B} \mathbf{A}^{-1}\mathbf{B}^\top$. Specifically, we assume that $\mathbf{A}$ and $\mathbf{S}$ are symmetric positive definite so that the preconditioner may be applied within the MINRES algorithm~\cite{paigesaunders}. If $\mathbf{C}$ is also symmetric positive semi-definite, we would have that the eigenvalues of the preconditioned matrix are contained within $[-1, \frac{1}{2}(1-\sqrt{5})] \cup [1, \frac{1}{2}(1+\sqrt{5})]$ regardless of problem dimension; see e.g., \cite[Cor. 1]{axelssonneytcheva2006}, \cite[Thm. 4]{pearsonthesis2013}, and \cite[Lem. 2.2]{silvesterwathen1994}. 
Now, due to the structure of $\bm{\mathcal{A}}$, applying $\mathbf{S}^{-1}$ can be very expensive due to $\mathbf{S}$ typically being a dense matrix. 
To offer a more favourable alternative to applying $\mathbf{S}^{-1}$ explicitly, one option is to approximate the Schur complement using a suitable product of matrices, the so-called ``matching strategy'' which was devised for simpler PDE-constrained optimization problems in~\cite{pearsonwathen2012,pearsonwathen2013}. We start by approximating the Schur complement with
\begin{equation*}
	\mathbf{S} \approx \mathbf{\hat{S}} := (\mathbf{B}+\mathbf{D}) \mathbf{A}^{-1} (\mathbf{B}+\mathbf{D})^\top,
\end{equation*}
where $\mathbf{D}$ is some real matrix to be determined.
Since we can rewrite the approximated Schur complement as 
$\mathbf{\hat{S}} = [ \left(\mathbf{B}+\mathbf{D}\right)\mathbf{A}^{-\frac{1}{2}} ] [ \left(\mathbf{B}+\mathbf{D}\right)\mathbf{A}^{-\frac{1}{2}} ]^\top$, we find that $\mathbf{\hat{S}}$ must be symmetric positive semi-definite, and provided $\mathbf{D}$ is suitably chosen (as below) $\mathbf{\hat{S}}$ will be symmetric positive definite. As this property also holds for the block $\mathbf{A}$, the preconditioner is then symmetric positive definite, which allows us to consider more options for iterative solvers shown later in this section. 
Expanding $\mathbf{\hat{S}}$ and matching the terms with those in $\mathbf{S}$, we determine $\mathbf{D}$ by specifying that  $\mathbf{D} \mathbf{A}^{-1} \mathbf{D}^\top$ should approximate $\mathbf{C}$. 
We can again rewrite this quantity as $[ \mathbf{D} \mathbf{A}^{-\frac{1}{2}} ] [ \mathbf{D} \mathbf{A}^{-\frac{1}{2}} ]^\top$, so the approximation of $\mathbf{C}$ must be a symmetric and positive semi-definite matrix.
A suitable approximation of $\mathbf{C}$ can be a block diagonal matrix, incorporating components that include mass matrices, so that
\begin{equation*}
	\mathbf{C} \approx \mathbf{\hat{C}} = 
	\text{blkdiag} \left( \left\lbrace \tau \alpha_1 \mathbf{M} \right\rbrace _{i=1}^{N_t-1}, \frac{\tau \alpha_1}{2} \mathbf{M}, 
	\left\lbrace \tau \alpha_2 \mathbf{M} \right\rbrace _{i=1}^{N_t-1}, \frac{\tau \alpha_2}{2} \mathbf{M} \right).
\end{equation*}
This is a convenient choice for $\mathbf{\hat{C}}$ because it captures many of the (linear) terms within the PDEs, and is also positive definite and hence approximable by a term of the form $\mathbf{D} \mathbf{A}^{-1} \mathbf{D}^\top$. Now we need to find the matrix $\mathbf{D}$. Since the block $\mathbf{A}$, and consequently, its inverse, is block diagonal, we obtain that $\mathbf{D}$ is a block diagonal matrix of the form
\begin{equation*}
	\mathbf{D}:=\text{blkdiag}\left(  \left\lbrace \mathbf{D}_{(1)}^i \right\rbrace_{i=1}^{N_t}  , \left\lbrace   \mathbf{D}_{(2)}^i\right\rbrace_{i=1}^{N_t} \right),
\end{equation*}
and we calculate that the diagonal blocks of $\mathbf{D}$ must be
\begin{equation*}
	\begin{aligned}
		\mathbf{D}_{(1)}^i &:= \tau \gamma \sqrt{\frac{\alpha_1}{\beta_1}} \mathbf{M} \quad \text{ for} \;  i=1,..., N_t-1, \hspace{5mm} 
		&\mathbf{D}_{(1)}^{N_t} &:= \tau \gamma \sqrt{\frac{\alpha_1}{2\beta_1}} \mathbf{M}, \\
		\mathbf{D}_{(2)}^i &:= \tau \gamma \sqrt{\frac{\alpha_2}{\beta_2}} \mathbf{M} \quad \text{ for} \;  i=1,..., N_t-1, 
		\hspace{5mm}     
		&\mathbf{D}_{(2)}^{N_t} &:= \tau \gamma \sqrt{\frac{\alpha_2}{2\beta_2}} \mathbf{M}.
	\end{aligned}
\end{equation*}

With an approximation of the ideal preconditioner $\bm{\hat{\mathcal{P}}}$, for computational efficiency we must define the solution of $\bm{\hat{\mathcal{P}}} \underline{\bm{y}} = \underline{\bm{z}}$ for a given $\underline{\bm{z}}$ without forming the inverse of the preconditioner explicitly. As $\bm{\hat{\mathcal{P}}}$ is block diagonal, we can split this computation into two parts. 
First, dealing with the $(1,1)$-block of $\hat{\bm{\mathcal{P}}}$ comes down to solving a system with scaled mass matrices on the diagonal. Instead of using a direct solver, we can use Chebyshev semi-iteration with Jacobi splitting, which provides an efficient approximation of mass matrices. For 2D problems discretized using P1 finite elements, as we later consider, we know that the eigenvalues of $\mathbf{M}$ preconditioned with its diagonal must lie between $0.5$ and $2$; see e.g., \cite{wathen1987}. We use this information within our implementation of Chebyshev semi-iteration and  let the method run for $20$ iterations.
Second, considering the approximated Schur complement within the $(2,2)$-block of the preconditioner, and with vector superscripts denoting the second halves of the vectors, we now aim to solve the system
$(\mathbf{B}+\mathbf{D}) \mathbf{A}^{-1} (\mathbf{B}+\mathbf{D})^\top \underline{\bm{y}}^2 = \underline{\bm{z}}^2$. As the matrix $\mathbf{B}+\mathbf{D}$ is block upper-triangular, we first solve $(\mathbf{B}+\mathbf{D}) \underline{\bm{m}} = \underline{\bm{z}}^2$ by block-backward substitution, after which we solve $ (\mathbf{B}+\mathbf{D})^\top \underline{\bm{y}}^2 = \mathbf{A} \underline{\bm{m}}$ by block-forward substitution. Each of the sub-systems can be solved using an algebraic multigrid solver (see \cite{vanek1996}, for instance), and we apply 6 V-cycles to ensure a sufficiently accurate approximation.

Finally, we highlight that there are a variety of preconditioned Krylov subspace iterative solvers for different classes of matrices and preconditioners. In particular, for symmetric matrices with symmetric and positive definite preconditioners of the structure \eqref{eq:ideal_prec}, the MINRES solver is the preferred choice. One could also apply the Bramble--Pasciak Conjugate Gradient Method \cite{bramblepasciak1988} with a modified preconditioner, and GMRES \cite{saadschultz}, BICG \cite{fletcher1976}, or its stabilized variant BiCGSTAB \cite{bicgstab} for problems involving non-symmetric systems or preconditioners, for example block triangular variants of \eqref{eq:ideal_prec}, which may also be readily applied. 
A key advantage of using such Krylov subspace solvers in general is their computational cost: they require multiplication with the system matrix, but this can be performed without having to form the matrix in full. That is, we are only required to define the linear operator that computes $\bm{\mathcal{A}}\underline{\bm{v}}$ for an arbitrary vector $\underline{\bm{v}}$ of consistent dimension. 
We may therefore utilize the specific block structure of the matrix under consideration, and define a matrix--vector product by simply using block-forward substitution for lower-triangular blocks and block-backward substitution for upper-triangular blocks.

\begin{remark}
To further motivate the above approximation of the Schur complement, consider the distribution of eigenvalues of $\mathbf{\hat{S}}^{-1} \mathbf{S}$ in a ``linearized'' setting, where $\bar{u} = \bar{v} = \bar{p} = \bar{q} = 0$. Then, $\mathbf{C}$ is positive definite, $\mathbf{D} \equiv \mathbf{C}^{\frac{1}{2}} \mathbf{A}^{\frac{1}{2}}$, $\mathbf{D} \mathbf{A}^{-1} \mathbf{D}^\top \equiv \mathbf{C}$, and
\begin{equation*}
    \mathbf{\hat{S}} = \mathbf{S} + \mathbf{C}^{\frac{1}{2}} \mathbf{A}^{-\frac{1}{2}} \mathbf{B}^\top + \mathbf{B} \mathbf{A}^{-\frac{1}{2}} \mathbf{C}^{\frac{1}{2}}.
\end{equation*}
By construction, both $\mathbf{S}$ and $\mathbf{\hat{S}}$ are symmetric positive definite. Hence, following the same reasoning as in \cite[Thm. 4.1]{pearsonwathen2013}, it may be shown that the eigenvalues of $\hat{S}^{-1} S$ are bounded below by $\frac{1}{2}$, irrespective of problem dimension and parameters involved. Further, it holds that
\begin{align*}
    \mathbf{C}^{\frac{1}{2}} \mathbf{A}^{-\frac{1}{2}} \mathbf{B}^\top + \mathbf{B} \mathbf{A}^{-\frac{1}{2}} \mathbf{C}^{\frac{1}{2}} ={}& \frac{\sqrt{\alpha_1 \beta_1}}{\gamma} \left( \begin{bmatrix}
	\mathbf{T}_1 & \mathbf{0} \\
	\mathbf{0} & \mathbf{0}
    \end{bmatrix} \otimes \mathbf{M} + \begin{bmatrix}
	\mathbf{T}_2 & \mathbf{0} \\
	\mathbf{0} & \mathbf{0}
    \end{bmatrix} \otimes \tau \mathbf{L}_{(1)} \right) \\
    & + \frac{\sqrt{\alpha_2 \beta_2}}{\gamma} \left( \begin{bmatrix}
	\mathbf{0} & \mathbf{0} \\
	\mathbf{0} & \mathbf{T}_1
    \end{bmatrix} \otimes \mathbf{M} + \begin{bmatrix}
        \mathbf{0} & \mathbf{0} \\
	\mathbf{0} & \mathbf{T}_2
    \end{bmatrix} \otimes \tau \mathbf{L}_{(2)} \right),
\end{align*}
where for the linearized problem $\mathbf{L}_{(1)} = \frac{D_u}{2} \mathbf{K} + \frac{\gamma}{2} \mathbf{M}$ and $\mathbf{L}_{(2)} = \frac{D_v}{2} \mathbf{K}$, and
\begin{equation*}
    \mathbf{T}_1 = \begin{bmatrix}
	2 & -1 \\
        -1 & 2 & \ddots \\
        & \ddots & \ddots & -1 \\
        & & -1 & 2 & -\frac{1}{\sqrt{2}} \\
        & & & -\frac{1}{\sqrt{2}} & \sqrt{2} \\
    \end{bmatrix}, \hspace{5mm} \mathbf{T}_2 = \begin{bmatrix}
	2 & 1 \\
        1 & 2 & \ddots \\
        & \ddots & \ddots & 1 \\
        & & 1 & 2 & \frac{1}{\sqrt{2}} \\
        & & & \frac{1}{\sqrt{2}} & \sqrt{2} \\
    \end{bmatrix}.
\end{equation*}
As $\mathbf{L}_{(1)}$, $\mathbf{L}_{(2)}$, $\mathbf{T}_1$, and $\mathbf{T}_2$ are all at least positive semi-definite, $\mathbf{C}^{\frac{1}{2}} \mathbf{A}^{-\frac{1}{2}} \mathbf{B}^\top + \mathbf{B} \mathbf{A}^{-\frac{1}{2}} \mathbf{C}^{\frac{1}{2}}$ is positive semi-definite, hence $\mathbf{S} \preceq \mathbf{\hat{S}}$ and the eigenvalues of $\mathbf{\hat{S}}^{-1} \mathbf{S}$ are bounded above by 1.
\end{remark}

\section{Comparison with backward Euler method}
\label{sec:backward_euler}
To evaluate the performance of our new numerical scheme, we compare it to the problem obtained by applying the backward Euler method for the time discretization of~\eqref{eq:OTD_Newt_st}--\eqref{eq:OTD_Newt_ad}; see \cite{stollpearsonmaini2016}.
For the Schnakenberg model, the linearized equations used in this section are \eqref{eq:SQP_state} and \eqref{eq:OTD_ad}. To reduce the number of variables, we can  substitute for the control variables using the gradient equations~\eqref{eq:OTD_gradient}.

Applying the backward Euler method over the time interval $(0,T)$, with $N_t$ time-steps of size $\tau$, we begin by discretizing the initial and final-time conditions,
\begin{equation}
		\mathbf{M}\bm{u}^0 = \bm{u}_0, \quad \mathbf{M}\bm{v}^0 = \bm{v}_0, \quad \mathbf{M}\bm{p}^{N_t} = \mathbf{M}\bm{q}^{N_t} = \bm{0},
	\label{eq:BWE_IC}
\end{equation}
with $\bm{u}_0$ and $ \bm{v}_0$ defined as in Section \ref{sec:otd_finite_elements}.
At each Newton iteration, discretization of the linearized PDEs leads to the following system:
\begin{equation}
	\begin{aligned}
		- \mathbf{M}\bm{u}^n + [\mathbf{M} + 2\tau \mathbf{L}_{(1)}^{n+1}] \bm{u}^{n+1} - \tau\gamma \mathbf{M}_{(u^{n+1})^2}\bm{v}^{n+1} - \frac{\tau \gamma^2}{\beta_1}\mathbf{M}\bm{p}^{n+1} = -\tau \bm{d}^{n+1}{}&, \\	
		2\tau\gamma \mathbf{M}_{u^{n+1} v^{n+1}} \bm{u}^{n+1} - \mathbf{M}\bm{v}^n  + [\mathbf{M} + 2\tau \mathbf{L}_{(2)}^{n+1}] \bm{v}^{n+1} - \frac{\tau \gamma^2}{\beta_2}\mathbf{M}\bm{q}^{n+1} = \tau \bm{d}^{n+1}{}&, \\
		[\tau \alpha_1 \mathbf{M} + \tau \mathbf{A}_{(1)}^{E,n}] \bm{u}^n +\tau \mathbf{A}_{(12)}^{E,n} \bm{v}^n + [\mathbf{M} + 2\tau \mathbf{L}_{(1)}^n]\bm{p}^n - \mathbf{M}\bm{p}^{n+1} + 2\tau\gamma \mathbf{M}_{u^n v^n} \bm{q}^n& \\
        = \tau \bm{c}^{n,n}{}&, \\
		\tau \mathbf{A}_{(12)}^{E,n}\bm{u}^n + \tau\alpha_2 \mathbf{M} \bm{v}^n - \tau\gamma \mathbf{M}_{(u^n)^2}\bm{p}^n +  [\mathbf{M} + 2\tau \mathbf{L}_{(2)}^n]\bm{q}^n - \mathbf{M}\bm{q}^{n+1} = \tau \bm{h}^{n,n}{}&, 
	\end{aligned}
	\label{eq:BWE_system}
\end{equation}
for $n=0,1,...,N_t-1$. 
Here we re-use the right-hand side vectors and some of the matrices introduced in Sections~\ref{sec:otd_finite_elements} and~\ref{sec:all_at_once}. The matrices specific to the backward Euler scheme are
\begin{equation*}
		\mathbf{A}_{(1)}^{E,i}:= 2 \gamma \mathbf{M}_{v^i(q^i - p^i)},  \quad 
		\mathbf{A}_{(12)}^{E,i}:=2 \gamma \mathbf{M}_{u^i(q^i - p^i)},  \quad \text{ for } \; 
		i=0,1,...,N_t-1.
\end{equation*}
As this system contains equations evolving forward and backward in time, we again solve for all the discretized variables at once until a stopping criterion for the Newton iterates is achieved. That is, we need to find solutions to the state variables $u$ and $v$ at time-steps $1,...,N_t$ and to the adjoint variables $p$ and $q$ at time-steps $0,1,...,N_t-1$.

We construct the system matrix for~\eqref{eq:BWE_system}, involving all four variables at time-steps $0,1,...,N_t$. Notice that
we can exclude the variables $\bm{p}^0$, $\bm{q}^0$, $\bm{u}^{N_t}$, and $\bm{v}^{N_t}$ from the all-at-once system, which can be done by decoupling the $N_t$-th equation of each of the state equations and the first equation of each of the adjoint equations. 
To obtain a symmetric system matrix, we reorder the variables, and using similar arguments as in Section~\ref{sec:precond}, we solve for the vector of unknowns $ \left[-\underline{\bm{p}},	\underline{\bm{u}} \right]$, where
\begin{equation*}
	\begin{aligned}
		\underline{\bm{u}} &:= \left[(\bm{u}^1)^{\top}, ..., (\bm{u}^{N_t-1})^{\top}, (\bm{v}^1)^{\top}, ..., (\bm{v}^{N_t-1})^{\top}\right]^{\top}, \\
		\underline{\bm{p}} &:= \left[(\bm{p}^1)^{\top}, ..., (\bm{p}^{N_t-1})^{\top}, (\bm{q}^1)^{\top}, ..., (\bm{q}^{N_t-1})^{\top}\right]^{\top},
	\end{aligned}
\end{equation*}
and with the system matrix
\begin{equation*}
	\bm{\mathcal{A}} := \\
	\begin{bmatrix}
\mathbf{A}^E
		& \begin{bmatrix}
			(\mathbf{B}_{11}^E)^\top & \mathbf{B}_{21}^E \\
			\mathbf{B}_{12}^E & (\mathbf{B}_{22}^E)^\top
		\end{bmatrix} \\
		\begin{bmatrix}
			\mathbf{B}_{11}^E & \mathbf{B}_{12}^E \\
			\mathbf{B}_{21}^E & \mathbf{B}_{22}^E
		\end{bmatrix}
		& \begin{bmatrix}
			-\mathbf{C}_{11}^E & -\mathbf{C}_{12}^E \\
			-\mathbf{C}_{12}^E & -\mathbf{C}_{22}^E 
		\end{bmatrix}
	\end{bmatrix}
	= 
	\begin{bmatrix}
		\mathbf{A}^E & \left(\mathbf{B}^E\right)^\top \\
		\mathbf{B}^E & -\mathbf{C}^E
	\end{bmatrix}	.
	\label{mat:BWE}
\end{equation*} 
The $(1,1)$-block and most of the sub-blocks are block diagonal:
\begin{equation*}
	\begin{aligned}
		\mathbf{A}^E &:= \text{blkdiag} \left( \hspace{-0.1 cm} \left\lbrace \frac{\tau \gamma ^2}{\beta_1} \mathbf{M}\right\rbrace_{i=1}^{N_t-1} ,	\left\lbrace \frac{\tau \gamma ^2}{\beta_2} \mathbf{M}\right\rbrace_{i=1}^{N_t-1} \right),  \\
		 \mathbf{C}_{11}^E &:= \text{blkdiag} \left(\hspace{-0.1 cm} \left\lbrace \tau \alpha_1 \mathbf{M} + \tau \mathbf{A}_{(1)}^{E,i}\right\rbrace _{i=1}^{N_t-1}\right),
   \end{aligned}
   \end{equation*}
   \vspace{-3mm}
   \begin{equation*}
   \begin{aligned}
        \mathbf{B}_{12}^E &:= \text{blkdiag} \left( \left\lbrace  2\tau\gamma	\mathbf{M}_{u^iv^i}\right\rbrace _{i=1}^{N_t-1}\right), 
	  &\mathbf{C}_{12}^E &:= \text{blkdiag} \left(  \left\lbrace \tau \mathbf{A}_{(12)}^{E,i}\right\rbrace _{i=1}^{N_t-1}\right),\\
		\mathbf{B}_{21}^E &:= \text{blkdiag} \left( \left\lbrace -\tau\gamma \mathbf{M}_{(u^i)^2}\right\rbrace _{i=1}^{N_t-1} \right),
	   &\mathbf{C}_{22}^E &:= \text{blkdiag} \left(  \left\lbrace \tau \alpha_2 \mathbf{M}\right\rbrace _{i=1}^{N_t-1}\right).	
	\end{aligned}
    \end{equation*}
Two of the sub-blocks appearing in the $(1,2)$- and $(2,1)$-blocks are upper-triangular: 
\begin{equation*}
	\begin{aligned}
		\mathbf{B}_{11}^E &:= \text{blkdiag} \left( \left\lbrace \mathbf{M} + 2\tau \mathbf{L}_{(1)}^i \right\rbrace _{i=1}^{N_t-1} \right) + \text{blksupdiag} \left(\left\lbrace -\mathbf{M} \right\rbrace_{i=1}^{N_t-2}  \right),  \\
		\mathbf{B}_{22}^E &:= \text{blkdiag} \left( \left\lbrace \mathbf{M} + 2\tau \mathbf{L}_{(2)}^i \right\rbrace _{i=1}^{N_t-1} \right) + \text{blksupdiag} \left(\left\lbrace -\mathbf{M} \right\rbrace_{i=1}^{N_t-2}  \right).
	\end{aligned}
\end{equation*}

Additionally, we need to solve the equations that were excluded from the system above. We solve the coupled system for $\bm{p}^0$ and $\bm{q}^0$:
\begin{flalign*}
\hspace{-1cm}
	&\begin{bmatrix}
		\mathbf{M} + 2\tau \mathbf{L}_{(1)}^0 &  2\tau\gamma \mathbf{M}_{u^0v^0} \\
		-\tau\gamma \mathbf{M}_{(u^0)^2}  & \mathbf{M} + 2\tau \mathbf{L}_{(2)}^0
	\end{bmatrix}
	\begin{bmatrix}
		\bm{p}^0 \\ \bm{q}^0
	\end{bmatrix}
        \hfill
        \\ 
        &= \begin{bmatrix}
		\tau\bm{c}^{0,0} + \mathbf{M}\bm{p}^1 - \tau[ \alpha_1 \mathbf{M} + \mathbf{A}_{(1)}^{E,0} ]  \bm{u}_0 - \tau \mathbf{A}_{(12)}^{E,0}\bm{v}_0 \\
		\tau\bm{h}^{0,0} + \mathbf{M}\bm{q}^1 - \tau \mathbf{A}_{(12)}^{E,0} \bm{u}_0 - \tau \alpha_2 \mathbf{M} \bm{v}_0
	\end{bmatrix},
 \end{flalign*}
and another system for $\bm{u}^{N_t}$ and $\bm{v}^{N_t}$:
\begin{equation*}
	\begin{bmatrix}
		\mathbf{M} + 2\tau \mathbf{L}_{(1)}^{N_t} & -\tau\gamma \mathbf{M}_{(u^{N_t})^2} \\
		 2\tau\gamma \mathbf{M}_{u^{N_t}v^{N_t}}    & \mathbf{M} + 2\tau \mathbf{L}_{(2)}^{N_t}  
	\end{bmatrix}
	\begin{bmatrix}
		\bm{u}^{N_t} \\ \bm{v}^{N_t}
	\end{bmatrix}
	=\begin{bmatrix}
		-\tau \bm{d}^{N_t} + \mathbf{M}\bm{u}^{N_t-1} \\
		\tau \bm{d}^{N_t} + \mathbf{M}\bm{v}^{N_t-1}
	\end{bmatrix},
\end{equation*}
using the vectors $\bm{p}^1$, $\bm{q}^1$, $\bm{u}^{N_t-1}$, and $\bm{v}^{N_t-1}$ previously calculated.

The matrix $\bm{\mathcal{A}}$ is again of (generalized) saddle-point form, and we proceed using a similar approach as in Section~\ref{sec:precond} and construct a block-diagonal preconditioner. 
In particular, the $(1,1)$-block of the preconditioner should approximate the $(1,1)$-block of the system matrix, which only has factors of mass matrices on the block diagonal. In the numerical implementation, we approximate the mass matrices using Chebyshev semi-iteration,  as in Section~\ref{sec:precond}. The $(2,2)$-block of the preconditioner is the approximation of the (negative) Schur complement $\mathbf{S}^E := \mathbf{C}^E + \mathbf{B}^E \left(\mathbf{A}^E \right)^{-1} \left( \mathbf{B}^E\right) ^\top$. 
The block $\mathbf{C}^E$ can be approximated by dropping the matrices emerging from nonlinear terms, 
leaving $\mathbf{\hat{C}}^E = \text{blkdiag} \left(\left\lbrace \tau \alpha_1 \mathbf{M} \right\rbrace_{i=1}^{N_t-1} \left\lbrace \tau \alpha_2 \mathbf{M} \right\rbrace_{i=1}^{N_t-1}  \right)$. 
We can thus approximate the Schur complement by $\mathbf{S}^E \approx \mathbf{\hat{S}}^E = (\mathbf{D}^E + \mathbf{B}^E) \left(\mathbf{A}^E \right)^{-1} \left( \mathbf{D}^E+ \mathbf{B}^E\right) ^\top$, and we find $\mathbf{D}^E$ by matching $\mathbf{D}^E \left(\mathbf{A}^E \right)^{-1} \left( \mathbf{D}^E\right) ^\top = \mathbf{\hat{C}}^E$ and neglecting the mixed terms after factorization of $\mathbf{\hat{S}}^E$, to obtain
\begin{equation*}
	\mathbf{D}^E := 
	\text{blkdiag} \left(\left\lbrace \tau \gamma \sqrt{\frac{\alpha_1}{\beta_1}} \mathbf{M}\right\rbrace _{i=1}^{N_t-1} \left\lbrace \tau \gamma \sqrt{\frac{\alpha_2}{\beta_2}}  \mathbf{M}\right\rbrace _{i=1}^{N_t-1} \right).
\end{equation*}
Since the system matrix is symmetric and positive definite, we use the preconditioned MINRES method to approximately solve the system. 
The equations for $\bm{p}^0$, $\bm{q}^0$, $\bm{u}^{N_t}$, and $\bm{v}^{N_t}$ can be solved with a direct solver, due to the relatively modest dimension of this system.

\section{Numerical results}
\label{sec:results}
To validate our numerical method, we now examine its performance on two relevant test problems. First, in Section~\ref{sec:convergence}, we benchmark our method and compare it with the backward Euler method described in Section~\ref{sec:backward_euler}, using a problem with an analytical solution. In Section~\ref{sec:datadriven}, we validate our approach for parameter identification in the context of a data-driven problem, with target states generated using the Schnakenberg equations with chosen functions $a$ and $b$. 

In all the following experiments we consider $\alpha := \alpha_1 = \alpha_2$ and,  without loss of generality, set $\alpha = 1$, and define $\beta:=\beta_1=\beta_2$. 
Further, we choose the diffusion parameters $D_u=1$, $D_v=10$. All the simulations are performed on a unit square domain $\Omega=(0,1)^2$. 
We consider zero-flux boundary conditions due to their more realistic physical significance, however, the Dirichlet boundary conditions can also be easily applied by adjusting the linear system.

To generate the matrices and the right-hand side vectors that emerge in the schemes, we use FEniCS \cite{LoggEtal2012} (version~2019.1.0), in particular the DOLFIN component. 
For other components of our implementation, we use Python \cite{python3manual} (v.~3.8.15) and its libraries: \verb*|scipy| \cite{scipy} (v.~1.6.3), \verb*|numpy| \cite{numpy} (v.~1.22.3), \verb*|matplotlib| \cite{matplotlib} (v.~3.4.0), and \verb*|pyamg| \cite{pyamg} (v.~4.1.0).

Using DOLFIN, we implement a finite element method in two dimensions with P$1$ elements and adhere to the corresponding node ordering. Once the matrices and vectors are generated, we convert them to Python arrays. The \verb*|scipy.sparse| library allows us to efficiently store the sparse matrices, while dense arrays and vectors are stored as \verb*|numpy| arrays. We also use the MINRES implementation from this library, which we provide with the system matrix and inverse of the preconditioner, both implemented as linear operators for memory efficiency. To perform the multigrid calculations in the Schur complement, we use the \verb*|smoothed_aggregation_solver| from the \verb*|pyamg| library. 
The plots in the following sections were generated with the \verb*|matplotlib.pyplot| library.

\subsection{Convergence scaling}
\label{sec:convergence}

\begin{table}[ht]
\begin{center}
	\caption{St\"{o}rmer--Verlet method for the Schnakenberg model: weighted errors, mean number of MINRES iterations and total SQP iterations, and CPU times, for different mesh sizes and $\beta=10^{-2}$.}
	\label{tab:SV_beta10-2}
	\begin{tabular}{ lllllllll }
		\toprule
		$i$ & DoF &$u_{\text{error}}$ &$v_{\text{error}}$ &$p_{\text{error}}$ &$q_{\text{error}}$ &$\bar{\text{it}}_M$ &$\text{it}_{\text{SQP}}$ &CPU \\
		\hline
		1 &24,200 &8.73e-2 &8.55e-2 &8.64e-3 &6.70e-3 &25 &6 &165 \\
		2 &176,400 &2.10e-2 &2.04e-2 &2.10e-3 &1.66e-3 &30 &4 &542 \\
		3 &1,344,800 &5.01e-3 &4.96e-3 &5.10e-4 &4.12e-4 &32 &4 &1,202 \\
		4 &10,497,600 &1.39e-3 &1.37e-3 &1.25e-4 &1.03e-4 &109 &4 &19,233 \\
		\bottomrule
	\end{tabular}
\end{center}
\end{table}

\begin{table}[ht]
\begin{center}
	\caption{St\"{o}rmer--Verlet method for the Schnakenberg model: weighted errors, mean number of MINRES iterations and total SQP iterations, and CPU times, for different mesh sizes and $\beta=10^{-3}$.}
	\label{tab:SV_beta10-3}
	\begin{tabular}{ lllllllll }
		\toprule
		$i$ & DoF &$u_{\text{error}}$ &$v_{\text{error}}$ &$p_{\text{error}}$ &$q_{\text{error}}$ &$\bar{\text{it}}_M$ &$\text{it}_{\text{SQP}}$ &CPU \\
		\hline
		1 &24,200 &4.61e-1 &2.04e-1 &6.18e-3 &2.92e-3 &30 &6 &169 \\
		2 &176,400 &1.69e-1 &7.24e-2 &1.73e-3 &8.44e-4 &38 &5 &410 \\
		3 &1,344,800 &5.29e-2 &2.25e-2 &4.66e-4 &2.37e-4 &38 &5 &1,745 \\
		4 &10,497,600 &1.39e-2 &5.99e-3 &1.15e-4 &6.06e-5 &74 &5 &16,605 \\
		\bottomrule
	\end{tabular}
\end{center}
\end{table}

\begin{table}[ht]
\begin{center}
	\caption{Backward Euler method for the Schnakenberg model: weighted errors, mean number of MINRES iterations and total SQP iterations, and CPU times, for different mesh sizes and $\beta=10^{-2}$.}
	\label{tab:BWE_beta10-2}
	\begin{tabular}{ lllllllll }
		\toprule
		$i$ & DoF &$u_{\text{error}}$ &$v_{\text{error}}$ &$p_{\text{error}}$ &$q_{\text{error}}$ &$\bar{\text{it}}_M$ &$\text{it}_{\text{SQP}}$ &CPU \\
		\hline
		1 &23,716 &1.03e-1 &9.53e-2 &8.13e-3 &6.90e-3 &30 &4 & 93\\
		2 &351,036 &2.47e-2 &2.25e-2 &1.96e-3 &1.74e-3 &30 &4 & 500\\
		3 &5,372,476 &5.92e-3 &5.47e-3 &4.77e-4 &4.30e-4 &25 &4 & 3,672\\
		4 &83,954,556 &1.23e-3 &1.29e-3 &1.13e-4 &9.90e-5 &27 &4 & 40,770\\
		\bottomrule
	\end{tabular}
\end{center}
\end{table}

\begin{table}[ht]
\begin{center}
	\caption{Backward Euler method for the Schnakenberg model: weighted errors, mean number of MINRES iterations and total SQP iterations, and CPU times, for different mesh sizes and $\beta=10^{-3}$.}
	\label{tab:BWE_beta10-3}
	\begin{tabular}{ lllllllll }
		\toprule
		$i$ & DoF &$u_{\text{error}}$ &$v_{\text{error}}$ &$p_{\text{error}}$ &$q_{\text{error}}$ &$\bar{\text{it}}_M$ &$\text{it}_{\text{SQP}}$ &CPU \\
		\hline
		1 &23,716 &5.09e-1 &2.14e-1 &5.95e-3 &2.38e-3 &29 &6 & 131\\
		2 &351,036 &1.98e-1 &7.97e-2 &1.76e-3 &8.40e-4 &36 &5 & 761\\
		3 &5,372,476 &6.49e-2 &2.61e-2 &4.96e-4 &2.50e-4 &36 &5 & 6,147\\
		4 &83,954,556 &1.95e-2 &7.78e-3 &1.39e-4 &7.25e-5 &34 &5 & 68,518\\
		\bottomrule
	\end{tabular}
\end{center}
\end{table}

To evaluate the performance of our method, we test it along with the backward Euler scheme on a problem with the following analytical solution:
\begin{equation*}
	\begin{aligned}
	u^\text{sol}(t,\textbf{x}) & = \text{e}^{0.1t} \left( \kappa(\textbf{x})+1\right), 
	&p^\text{sol}(t,\textbf{x}) &= \left( \text{e}^{0.1t} -\text{e}^{0.1T}  \right) \left( \kappa(\textbf{x}) +1 \right),  \\
	v^\text{sol}(t,\textbf{x}) & = \text{e}^{0.15t} \left( \eta(\textbf{x}) +1 \right), 
	\hspace{5mm} 
	&q^\text{sol}(t,\textbf{x}) &= \left( \text{e}^{0.15t} -\text{e}^{0.15T}  \right) \left( \eta(\textbf{x}) +1 \right). 
	\end{aligned}
\end{equation*}
Here, $\kappa(\textbf{x}) = \cos(2 \pi x_1) \cos(2 \pi x_2)$ and $\eta(\textbf{x}) = \cos(\pi x_1) \cos(\pi x_2)$, with $\textbf{x}=(x_1, x_2)$. The initial conditions for the state variables are constructed using the analytical solutions at $t=0$. For the desired states, we define the functions
\begin{equation*}
	\begin{aligned}
		\hat{u}^\text{sol}(t,\textbf{x}) ={}& \frac{1}{\alpha}\left[ -0.1\text{e}^{0.1t}(\kappa(\textbf{x}) + 1) + 		8D_u\pi^2\left( \text{e}^{0.1t} -\text{e}^{0.1T}  \right)\kappa(\textbf{x})
		+ \alpha u^\text{sol}(t,\textbf{x})  \right.\\
		&\left. {}+
		2\gamma u^\text{sol}(t,\textbf{x}) v^\text{sol}(t,\textbf{x}) (q^\text{sol}(t,\textbf{x})-p^\text{sol}(t,\textbf{x})) +
		\gamma p^\text{sol}(t,\textbf{x})) \right] , \\
		\hat{v}^\text{sol}(t,\textbf{x}) ={}& \frac{1}{\alpha}\left[ -0.15\text{e}^{0.15t}(\eta(\textbf{x}) + 1) 
            + 2D_v\pi^2\left( \text{e}^{0.15t} -\text{e}^{0.15T}  \right)\eta(\textbf{x}) \right.\\
		&\left. {}+ \alpha v^\text{sol}(t,\textbf{x})  
		+ \gamma u^\text{sol}(t,\textbf{x})^2 (q^\text{sol}(t,\textbf{x})-p^\text{sol}(t,\textbf{x})) \right].
	\end{aligned}
\end{equation*}
Finally, we need to introduce source functions $f^\text{sol}(t,\textbf{x})$ and $g^\text{sol}(t,\textbf{x})$ to the right-hand side of the first and second state equations in \eqref{eq:schnakenberg_eq}, respectively:
\begin{equation*}
	\begin{aligned}
		f^\text{sol}(t,\textbf{x}) ={}& (0.1 + \gamma)u^\text{sol}(t,\textbf{x}) + 8D_u\pi^2 \text{e}^{0.1t}\kappa(\textbf{x}) - \gamma |u^\text{sol}(t,\textbf{x})|^2 v^\text{sol}(t,\textbf{x}) \\ 
        &{}- \frac{\gamma^2}{\beta} p^\text{sol}(t,\textbf{x}), \\
		g^\text{sol}(t,\textbf{x}) ={}& 0.15v^\text{sol}(t,\textbf{x}) + 2D_v\pi^2 \text{e}^{0.15t}\eta(\textbf{x}) + \gamma |u^\text{sol}(t,\textbf{x})|^2 v^\text{sol}(t,\textbf{x}) 
        \\ &{}- \frac{\gamma^2}{\beta} q^\text{sol}(t,\textbf{x}).
	\end{aligned}
\end{equation*}
In the notation of \eqref{eq:schnakenberg_eq}, this corresponds to taking $\Phi(u,v) = \gamma(u-u^2 v) - f^\text{sol}$ and $\Psi(u,v) = \gamma u^2 v - g^\text{sol}$, which does not alter the adjoint or gradient equations derived thereafter.

As both schemes yield a symmetric system matrix, we use the MINRES solver. The preconditioner is built using the approach described in Sections~\ref{sec:precond} and~\ref{sec:backward_euler}. The MINRES solver terminates when the relative residual is below $\text{tol}_\text{M}$. We found that for our method, $\text{tol}_\text{M} = 10^{-8}$ is sufficient for $\beta=10^{-3}$, however, we report results with $\text{tol}_\text{M} = 10^{-9}$ as this tolerance was required to achieve the correct scaling of the error for $\beta=10^{-2}$. The outer iteration, i.e., the SQP method for the St\"{o}rmer--Verlet scheme or Newton method for the backward Euler scheme, stops after the relative error with the solution from the previous iteration reaches $\text{tol}_{\text{SQP}}=10^{-5}$ in all variables.

The first iteration of the SQP/Newton method is commenced using suitable initial guesses. In the case of the St\"{o}rmer--Verlet scheme, we first solve the problem on the coarsest mesh where our initial guesses for the state variables are taken to be the desired states, while we take vectors of zeros as the initial guesses for the adjoint variables. After obtaining the solution, we proceed with solving on a finer mesh, using the previously calculated solution interpolated onto the current mesh and time interval and multiplied by $0.8$ as an initial guess. This multiplication is to ensure that we are able to retrieve an accurate solution, without initializing in the immediate neighbourhood of the solution. We repeat this process for finer meshes. To accommodate the inbuilt MINRES solver, we use slightly different initial guesses for the backward Euler scheme on the coarsest grid level: specifically, we use the target states multiplied by $0.4$ as initial guesses for $u$, $v$. For finer meshes, we proceed in the same way as outlined for the St\"{o}rmer--Verlet scheme. We note that using this approach for the St\"{o}rmer--Verlet scheme leads to the same numerical solution.
 
The results in Tables \ref{tab:SV_beta10-2}--\ref{tab:BWE_beta10-3} were computed on a CentOS Linux HPC facility with Intel(R) Xeon(R) Gold 6348 CPU @ 2.60GHz.
The meshes are defined by the spacing $h=\frac{2^{1-i}}{10}$ for $i\in\{1,2,3,4\}$, with $\tau=\frac{h}{5}$ for the St\"{o}rmer--Verlet scheme and $\tau=2h^2$ for the backward Euler scheme, with the total number of time-steps $N_t^S$ and $N_t^E$, respectively.
The number of degrees of freedom (DoF) is given by $4N_t^SN_x$ for the St\"{o}rmer--Verlet and $4(N_t^E-1)N_x$ for the backward Euler scheme.
We state the weighted errors with the exact solution for $T=1$, $\gamma=2$, $\beta \in \{10^{-2},10^{-3}\}$. 
We report the $\ell^2$  in space and $\ell^\infty$ in time error, given by 
\begin{equation*}
	u_{\text{error}} := \max_{ 0 \leq i \leq N_t} (h |\bm{u}^i - u^\text{sol}(i\tau,\textbf{x}) |_2 ),
\end{equation*}
{\c using values of the exact solution on nodal points}, which is used for all the variables. 
The CPU time measurements only include the SQP/Newton method and the calculation of the system that is solved with MINRES, i.e., we do not include the operations preceding the linear solve nor the calculation of the variables that were excluded from the large system, i.e., $\bm{p}^0$, $\bm{q}^0$ and $\bm{u}^{N_t}$, $\bm{v}^{N_t}$, in the time measurement. 

Comparing the accuracy between the two schemes in Tables \ref{tab:SV_beta10-2}--\ref{tab:BWE_beta10-3}, we notice that for the same mesh size, both schemes achieve errors of roughly the same order. Comparing the errors across the two values of $\beta$ considered, our strategy generally outperforms the backward Euler method. 
Note in particular that for $\beta=10^{-3}$ at the finest mesh, the St\"ormer--Verlet scheme achieves smaller errors than the backward Euler scheme despite having eight times larger time-step size.

While the number of SQP iterations in Tables \ref{tab:SV_beta10-2} and \ref{tab:SV_beta10-3} is fairly consistent across the different meshes, the number of MINRES iterations increases about three-fold for $\beta=10^{-2~}$ and two-fold for $\beta=10^{-3}$ at the finest mesh level. 
This increase might be due to the denser structure of the system matrix resulting from our strategy applied to the PDECO problem, as opposed to the backward Euler scheme, combined with the large dimension of the problem and the features of the nonlinearity of the problem affecting the linear algebraic properties of the system at this mesh level. Despite this, we can noticeably reduce the CPU time for both values of $\beta$, while achieving similar accuracy as with the backward Euler scheme. Therefore, using our strategy for solving the given PDECO problem at a fine resolution allows us to save computational resources since we can reduce the problem dimension by taking larger time-steps while maintaining good accuracy.

\subsection{Data-driven problem}
\label{sec:datadriven}

\begin{table}[ht]
\begin{center}
	\caption{Parameter identification problem: squared $L^2(Q)$ norms of individual terms in cost functional, average number of MINRES iterations, and total SQP iterations, with $\beta=10^{-2}$.}
		\begin{tabular}{llllllll}
					\toprule
			$i$     & $L^2(Q) [u-\hat{u}]^2$ & $L^2(Q) [v-\hat{v}]^2 $ & $L^2(Q) [a] ^2$ & $L^2(Q) [b] ^2$ & $\bar{\text{it}}_M$ &$\text{it}_{\text{SQP}}$ \\
						\hline
			2   & 2.89e-5               & 5.74e-5                                                          & 0.1786          & 0.1911          & 14          & 4          \\
			3  & 2.85e-5               & 4.16e-5                                                          & 0.1789          & 0.1877          & 13           & 4          \\
			4 & 2.89e-5               & 3.75e-5                                                          & 0.1791          & 0.1890          & 12          & 8         \\
			\hline
		\end{tabular}
		\label{tab:dd_beta10-2}
\end{center}
\end{table}

\begin{table}[ht]
\begin{center}
	\caption{Parameter identification problem: squared $L^2(Q)$ norms of individual terms in cost functional, average number of MINRES iterations, and total SQP iterations, with $\beta=10^{-3}$.}
	\begin{tabular}{llllllll}
		\toprule
		$i$    & $L^2(Q) [u-\hat{u}]^2$ & $L^2(Q) [v-\hat{v}]^2 $ & $L^2(Q) [a] ^2$ & $L^2(Q) [b] ^2$ & $\bar{\text{it}}_M$ &$\text{it}_{\text{SQP}}$  \\
		\hline
		2   & 3.03e-7               & 2.25e-6                                                          & 0.1808          & 0.2057          & 11          & 4          \\
		3  & 3.31e-7               & 1.04e-6                                                          & 0.1810          & 0.1990          & 10          & 4          \\
		4 & 4.72e-7               & 7.96e-7                                                          & 0.1806          & 0.1997          & 9          & 6         \\
		\hline
	\end{tabular}
	\label{tab:dd_beta10-3}
\end{center}
\end{table}

\begin{table}[ht]
\begin{center}
		\caption{Parameter identification problem: squared $L^2(Q)$ norms of individual terms in cost functional, average number of MINRES iterations, and total SQP iterations, with $\beta=10^{-4}$.}
	\begin{tabular}{llllllll}
		\toprule
		$i$     & $L^2(Q) [u-\hat{u}]^2$ & $L^2(Q) [v-\hat{v}]^2 $ & $L^2(Q) [a] ^2$ & $L^2(Q) [b] ^2$ & $\bar{\text{it}}_M$ &$\text{it}_{\text{SQP}}$  \\
		\hline
		2   & 3.42e-9               & 8.84e-8                                                          & 0.1811          & 0.2110          & 8          & 4           \\
		3  & 4.35e-9               & 9.00e-8                                                          & 0.1812          & 0.2015          & 8          & 4           \\
		4 & 6.30e-8               & 8.39e-8                                                          & 0.1813          & 0.2014          & 6          & 4          \\
		\hline
	\end{tabular}
	\label{tab:dd_beta10-4}
\end{center}
\end{table}

\begin{figure}[ht]
\begin{center}
	\caption{Snapshot of the solution at $t=0.7$ for the desired states $\hat{u}$, $\hat{v}$ (left) and computed state variables $u$, $v$ (middle), and at $t=0.695$ for the computed control variables $a$, $b$ (right), with $\beta=10^{-2}$. The colorbar takes account of values from the whole time interval.}
	\includegraphics[width=\textwidth]{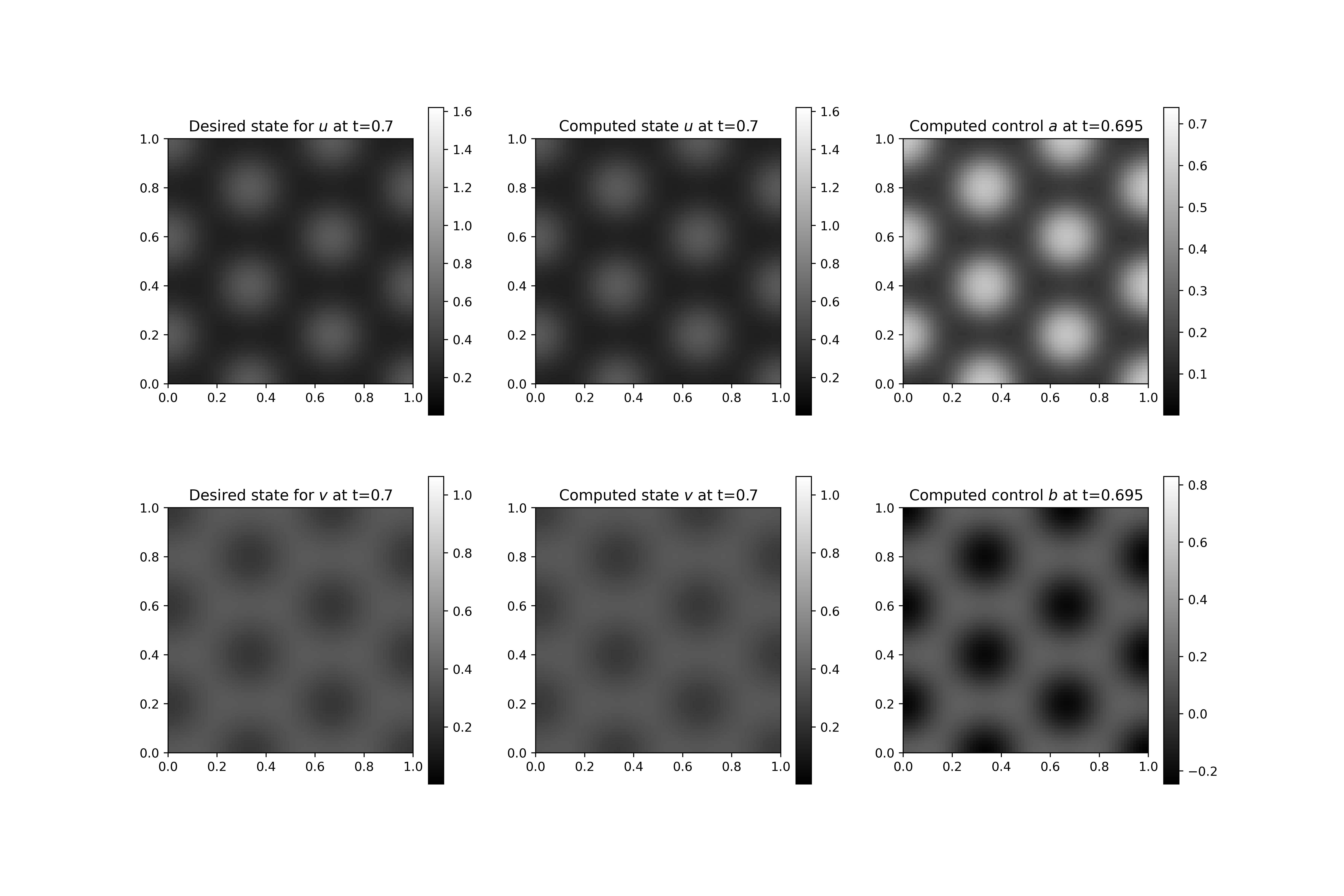}
	\label{fig:dd_t0.7}
\end{center}
\end{figure}

\begin{figure}[ht]
\begin{center}
	\caption{Snapshot of the solution at $t=2$ for the desired states $\hat{u}$, $\hat{v}$ (left) and computed state variables $u$, $v$ (middle), and at $t=1.995$ for the computed control variables $a$, $b$ (right), with $\beta=10^{-2}$. The colorbar takes account of values from the whole time interval.}
	\label{fig:dd_t2}
	\includegraphics[width=\textwidth]{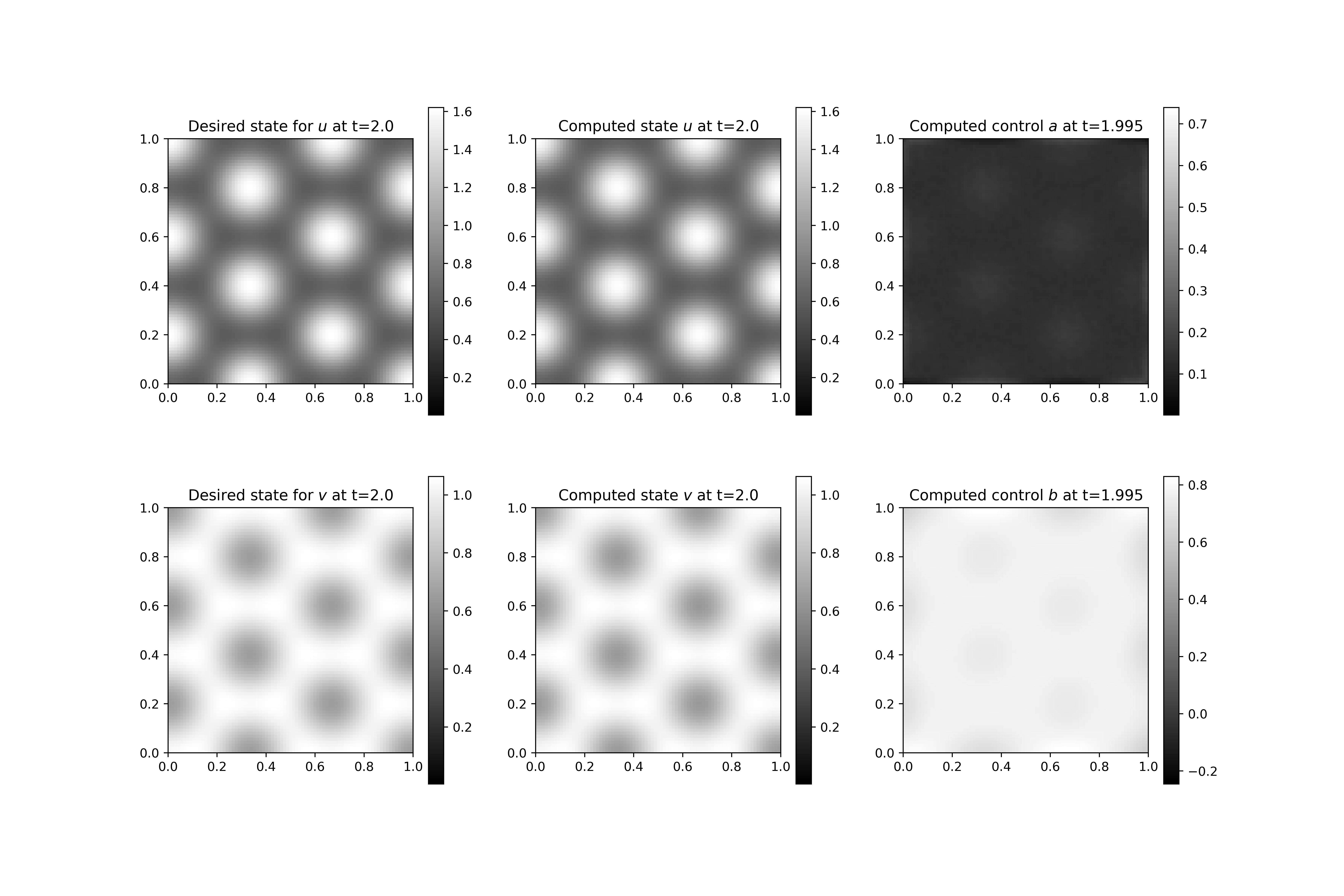}
\end{center}
\end{figure}

\begin{figure}[ht]
\begin{center}
	\centering
	\caption{Evolution of the mean of the numerical solution for the control variables $a$, $b$ across the time interval compared to those of the control variables $a_G$, $b_G$ used to generate the target state at $t=5$, with $\beta=10^{-2}$.}
	\label{fig:dd_means}
	\includegraphics[scale=0.5]{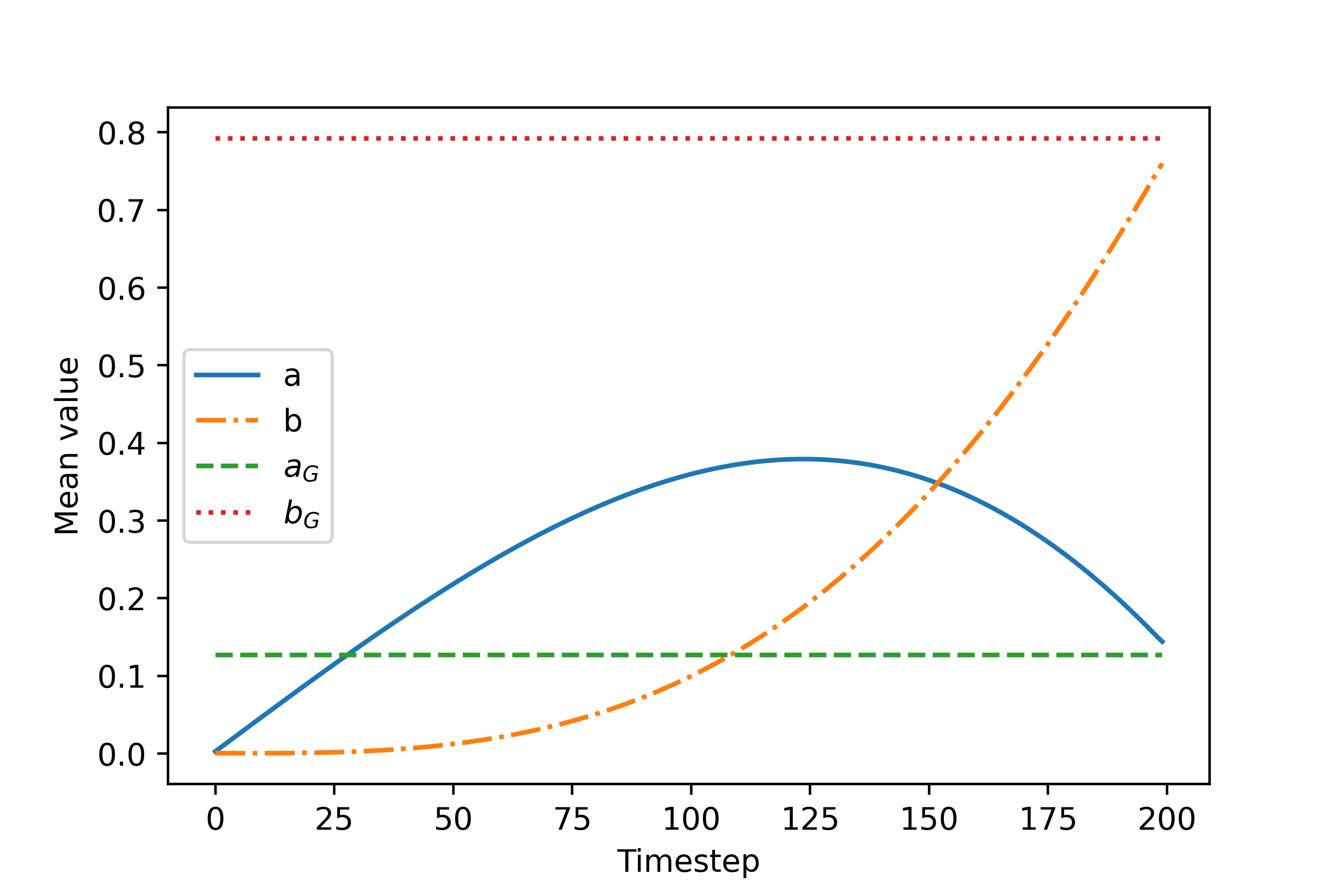}
\end{center}
\end{figure}

For testing purposes, we generate the target states $\hat{u}$, $\hat{v}$ directly using the Schnakenberg model implemented in FEniCS.
For the Turing patterns to emerge, we introduce small perturbations to the homogeneous steady states as in \cite{garvie2010}, with values for the control parameters $a_G=0.126779$ and $b_G = 0.792366$ taken to be constant across the whole domain. We take snapshots of the patterns once they are stationary, at time $t=5$. 
Since in our PDECO problem we assume the target states evolve in time, for simplicity we linearly interpolate from target states that are zero at time $t=0$ to acquire the snapshots at time $t=T$.
Here we used $T=2$, the fixed time-step size $\tau = 10^{-2}$, and tolerances $\text{tol}_{\text{M}} = 10^{-7}$ and $\text{tol}_{\text{SQP}} = 10^{-6}$.

In Tables \ref{tab:dd_beta10-2}--\ref{tab:dd_beta10-4}, we observe the squared $L^2(Q)$-norms which form the cost functional of the PDECO problem along with the iteration numbers, as in the previous section. 
The squared $L^2(Q)$-norms of the misfits between the desired and computed states remain very small and similar for different mesh sizes across different values of $\beta$. The norms of the control parameters also remain similar in each table, which suggests consistency of the results regardless of the mesh resolution.
The average number of MINRES iterations decreases for finer meshes; however, as the size of the problem increases, more SQP iterations are needed, which differs from the behaviour observed in the previous section when the time-step size was varied. 

The outcome of the numerical simulations can be observed in Figures~\ref{fig:dd_t0.7} and~\ref{fig:dd_t2} with snapshots of the desired states, computed states and computed controls at $t=0.7$ and the final time $t=2$, respectively. Both the desired and computed states in each snapshot are in good agreement which is confirmed by the norms of the differences $u - \hat u$, $v - \hat v$ in Table~\ref{tab:dd_beta10-2}. We note that since no precautions were taken for positivity preservation, the control variable $b$ has reached negative values within the time interval. 
Importantly, it is obvious that to drive the state variables to the desired pattern, the values of the control variables cannot be homogeneous across the spatial domain. 
Therefore we compare the mean values of the control variables across the domain with the constants $a_G$ and $b_G$, used to generate the images at the final time. In Figure \ref{fig:dd_means}, we observe that at the final time the means of the computed controls do indeed roughly achieve the desired control values. Since the desired state to be achieved at final time $T=2$ was generated using the state equations at time $t=5$, we can also conclude that we were able to drive the state to the desired state in a shorter time. 
This means that even though in reality we may not know the actual duration of the evolution into a particular pattern for a biological process, we can still obtain an accurate approximation for the unknown parameters.

\section{Conclusions}\label{sec:conc}

In this paper, we have devised an SQP-based method for challenging and large-scale parameter identification problems constrained by reaction--diffusion equations, based on the Schakenberg model for pattern formation in mathematical biology. Specifically, we have devised a paradigm for time-stepping and approximation of cost functional terms, used in conjunction with a suitable resolution in the spatial variables, such that optimization and discretization steps commute within our numerical algorithm. Further, our method demonstrates second-order convergence upon mesh refinement in time and space. We devise feasible and efficient preconditioned iterative methods for solving the linear systems arising from each iteration of the SQP solver. Problems with dimensions of the order of tens of millions were successfully tackled using our solution strategy. Future work includes devising a theoretical stability and convergence analysis for our numerical scheme and applying it to even more demanding problems from mathematical biology.

\noindent \textbf{Acknowledgments.}
This work has made use of the resources provided by the Edinburgh Compute and Data Facility (ECDF) (\url{http://www.ecdf.ed.ac.uk}).

\noindent \textbf{Funding.}
KB was supported by the EPSRC Centre for Doctoral Training in Mathematical 
Modelling, Analysis and Computation (MAC-MIGS) funded by the UK Engineering and 
Physical Sciences Research Council (EPSRC) grant EP/S023291/1, Heriot-Watt University, and the 
University of Edinburgh. JWP was supported by the EPSRC grant EP/S027785/1.

\noindent \textbf{Code availability.}
The code is publicly available on GitHub at \url{https://github.com/KarolinaBenkova/PDECO-schnak}.

\noindent \textbf{Data availability.}
The data generated for the data-driven problem can be found on GitHub along with the code.

\noindent \textbf{Conflict of interest.} The authors declare that they have no conflict of interest.

\bibliographystyle{ieeetr}
\bibliography{literature}

\end{document}